\documentclass[11t]{article}
\usepackage[T1]{fontenc}
\usepackage[utf8]{inputenc}
\usepackage[cyr]{aeguill}
\usepackage[english]{babel}
\usepackage{amsmath,amssymb,amsthm}
\usepackage{bm}
\usepackage{graphicx}

\usepackage{authblk}	
\usepackage{geometry}
\usepackage{mathrsfs}

\usepackage{multirow}
\usepackage[justification=centering]{caption}
\usepackage{subfig}
\usepackage{xcolor,colortbl}
\usepackage{diffcoeff}

\usepackage[colorlinks=true,linkcolor=black, citecolor=blue, urlcolor=blue]{hyperref}

\DeclareMathOperator{\Tr}{Tr}

\newtheorem{theorem}{Theorem}
\newtheorem{remark}{Remark}
\newtheorem{proposition}{Proposition}
\newtheorem{definition}{Definition}
\newtheorem{corollary}{Corollary}

\newcommand{\sqboxs}{1.2ex}
\newcommand{\sqbox}[1]{\textcolor{#1}{\rule{\sqboxs}{\sqboxs}}}

\definecolor{Blue}{rgb}{0,1,1}
\definecolor{Green}{rgb}{0.4,1,0.4}
\definecolor{GreenYellow}{rgb}{0.3,0.6,0.5}
\definecolor{Yellow}{rgb}{1,1,0.4}
\definecolor{Orange}{rgb}{1,0.6,0.2}
\definecolor{Purple}{rgb}{1,0.5,0.5}
\definecolor{Red}{rgb}{1,0.2,0.2}

\definecolor{Grey}{rgb}{0.5,0.5,0.5}

\newcolumntype{g}{>{\columncolor{Green}[.95\tabcolsep]}c}
\newcolumntype{m}{>{\columncolor{Yellow}[.95\tabcolsep]}c}
\newcolumntype{b}{>{\columncolor{Blue}[.95\tabcolsep]}c}

\graphicspath{{./Figures/}, {./Figures/Eig_Mindlin/}, {./Figures/Simulations/}}

\def\onedot{$\mathsurround0pt\ldotp$}
\def\cddot{
	\mathbin{\vcenter{\baselineskip.67ex
			\hbox{\onedot}\hbox{\onedot}}%
}}

\makeatletter \renewcommand\d[1]{\ensuremath{%
		\;\mathrm{d}#1\@ifnextchar\d{\!}{}}}
\makeatother

\title{Port-Hamiltonian formulation and \\ Symplectic discretization of plate models\\
	Part I : Mindlin model for thick plates}	

\author[1]{Andrea Brugnoli\thanks{andrea.brugnoli@isae.fr}}
\author[1]{Daniel Alazard\thanks{daniel.alazard@isae.fr}}
\author[1]{Val\'erie Pommier-Budinger \thanks{valerie.budinger@isae.fr}}
\author[1]{Denis Matignon \thanks{denis.matignon@isae.fr}}
\affil[1]{ISAE-SUPAERO, Universit\'e de Toulouse, France. \\
	10 Avenue Edouard Belin, BP-54032, 31055 Toulouse Cedex 4.}

\begin{document}
\maketitle
	
	\begin{abstract}
			The port-Hamiltonian formulation is a powerful method  for modeling and interconnecting systems of different natures. In this paper, the port-Hamiltonian formulation in tensorial form of a thick plate described  by the Mindlin-Reissner model is presented. Boundary control and observation are taken into account. Thanks to tensorial calculus, it can be seen that the Mindlin plate model mimics the interconnection structure of its one-dimensional counterpart, i.e. the Timoshenko beam.\\
			The Partitioned Finite Element Method (PFEM\footnote{PFEM stands for partitioned finite element method.}) is then extended to both the vectorial and tensorial formulations in order to obtain a suitable, i.e. structure-preserving, finite-dimensional port-Hamiltonian system (PHs\footnote{PHs stands for port-Hamiltonian systems.}), which preserves the structure and properties of the original distributed parameter system. Mixed boundary conditions are finally handled by introducing some algebraic constraints. \\
			Numerical examples are finally presented to validate this approach. 
		\end{abstract}

	\section*{Introduction}
	The Hamiltonian formalism arising from the Legendre transformation of Euler-Lagrange system  has been already widely explored \cite{SymplecticElasticity}. The Legendre transformation gives rise to a system of equations ruled by a Hamiltonian matrix on a symplectic space. Port-Hamiltonian systems (PH) are instead defined with respect to a {Hamiltonian operator \cite[Chapter~7]{Olver}. For a complete discussion on as the relation between Lagrangian and Hamiltonian formulation the reader may consult \cite{MardsenDiraclAG_I, MardsenDiraclAG_II}.} The PH framework is acquiring more visibility for its capability to represent a huge class of systems coming from different realms of physics in a modular way. Finite-dimensional port-Hamiltonian systems can be easily interconnected together, as shown in \cite{Cervera2007}, allowing the construction of complex multi-physics systems. The interconnection is possible also in the infinite-dimensional case \cite{ShaftIntInfinite}, even if the procedure is not as straightforward as in the finite-dimensional case. Eventually, it is also possible to merge finite and infinite PH systems \cite{vanderShaftintFinInf}. These features and capabilities are particularly appealing for control engineers in order to simplify the modeling task in preliminary analyses.  \\

	Distributed parameter systems are of relevant interest given the increased computational power available for simulations. PH distributed systems were initially presented in \cite{VANDERSCHAFT2002166}, by using the theory of differential forms. Links towards functional analysis have been made in \cite{Villegas} and an exhaustive reference about the subject can be found in \cite{BookZwart}. The fundamental feature of a distributed PH system is the underlying geometric interconnection structure, the so-called Stokes-Dirac structure, that describes the power flow across the boundary and inside the system, together with an energy functional, the Hamiltonian, that determines the nature of the system. Linear/nonlinear, parabolic/hyperbolic systems can be all recast into this framework, \cite{bookPHs}. Port-Hamiltonian systems are by definition open systems, able to interact with the environment through boundary ports. The definitions of these boundary variables is of utmost importance to show that a PH system is well posed (see \cite{LeGorrec2005} for the proof in the 1D case). Applications coming from continuum mechanics, electrodynamics and thermodynamics can be integrated inside this framework. Academic examples typically considered are the transmission line, the shallow water equations and the Maxwell equations \cite{VANDERSCHAFT2002166}. \\
	
	Numerical simulations and control techniques require a spatial discretization that is meant to preserve the underlying properties related to power continuity. The discretization procedure for systems under PH formalism consists of two steps:
	\begin{itemize}
		\item Finite-dimensional approximation of the Stokes-Dirac structure, i.e. the formally skew symmetric differential operator that defines the structure. The duality of the power variables has to be mapped onto the finite approximation. The subspace of the discrete variables will be represented by a Dirac structure. 
		\item The Hamiltonian requires as well a suitable discretization, which gives rise to a discrete Hamiltonian. 
	\end{itemize} 
	
	The research community is focusing on structure-preserving discretization techniques since several years. In \cite{Golo}, the authors made use of a mixed finite element spatial discretization for 1D hyperbolic system of conservation laws, introducing different low-order basis functions for the energy and co-energy variables. Pseudo-spectral methods relying on higher-order global polynomial approximations were studied in \cite{moulla:hal-01625008}. This method was used and extended to take into account unbounded control operators in \cite{articleFlavio}. More recently a simplicial discretization based on discrete exterior calculus was proposed in \cite{SESLIJA20121509}. This approach comes with additional complexities, since a primal and a dual meshes have to be defined but the discretization is structure-preserving, regardless of the spatial dimension of the problem. Weak formulations which lead to Galerkin numerical approximations began to be explored in the last years. In \cite{WeakForm_Kot} the prototypical example of hyperbolic systems of two conservation law was discretized by a weak formulation. In this approach the boundary is split according to the causality of boundary ports, so that mixed boundary conditions are easily handled, but still dual meshes have to be defined.  \\

	{The symplectic approach based on the Legendre transformation has been used to deal with the Mindlin plate model in \cite{WeianMindlin}. It has been capable of providing insights on analytical solutions for the eigenproblem of rectangular plates. This methodology is of use whenever easy engineering solutions are sought after. On the contrary, the port-Hamiltonian framework is a powerful tool whenever complex systems have to be modeled in a structured way.} The main contribution of this paper is the enrichment of the PH formulation for the Mindlin plate model, by making use of tensor calculus \cite[Chapter~16]{Grinfield}. This kind of model was already presented in \cite{MacchelliMindlin} and {using the jet bundle formalism in \cite{jetMin}} but here a novel formulation based on tensorial calculus is introduced.  The second contribution of this paper concerns the extension of the Partitioned Finite Element Method (PFEM) to the Mindlin plate model. In this approach, originally presented in \cite{CardosoRibeiro2018}, once the system has been put into weak form, a subset of the equations is integrated by parts, so that boundary variables are naturally included into the formulation and appear as control inputs, the collocated outputs being defined accordingly. Then, the discretization of energy and co-energy variables (and the associated test functions) leads directly to a full rank representation for the finite-dimensional port-Hamiltonian system. If mixed boundary conditions are to be considered, the finite-dimensional system contains constraints, leading to an algebraic differential system (DAE), which can be treated and analyzed by referring to \cite{vanderSchaft2013, beattie2018linear}. This approach, similarly to the one detailed in \cite{WeakForm_Kot}, makes possible the usage of FEM software, like FEniCS \cite{LoggMardalEtAl2012}. The final discretized system can be further reduced by using appropriate model reduction techniques, as explained in \cite{TONG20132727, Mehrmann2018}. \\
	
	The paper is divided into {five} main sections 
	\begin{enumerate}
		\item The framework of finite-dimensional port-Hamiltonian systems (PHs) and the notion of Dirac and Stokes-Dirac structures are recalled. The infinite-dimensional case is then illustrated by means of the Timoshenko beam model, the 1D counterpart of the Mindlin plate model.
		\item {The constitutive relations, equations of motions and energies for the Mindlin-Reissner plate are detailed.}
		\item The Mindlin plate port-Hamiltonian formulation is highlighted by defining energy and co-energy variables and the interconnection structure. The boundary variables are found by introducing the energy balance. Then the underlying Stokes-Dirac structure is easily obtained by defining the flow and effort spaces, together with the space of boundary variables.
		\item The PFEM discretization for the Mindlin plate is detailed. The problem is put into weak form first, then the necessary integrations by parts are performed and finally the basis functions for the energy, co-energy and test functions are chosen.
		\item {To demonstrate the consistency of the approach an eigenvalue study of a square plate considering various boundary conditions is performed, together with some time-domain simulation.} 
	\end{enumerate}

	{
		\section{Reminder on port-Hamiltonian systems}
		In this section the concepts of Dirac structure and Stokes-Dirac structure are recalled. From these concepts finite and infinite dimensional port-Hamiltonian system are then derived by considering the port behavior. 
		\subsection{Dirac Structures}
		Consider an $n$-dimensional space $\mathcal{F}$ over the field $\mathbb{R}$ and $\mathcal{E} \equiv \mathcal{F}'$ its dual, i.e. the space of linear operator $\bm{e} : \mathcal{F} \rightarrow \mathbb{R}$. The elements of $\mathcal{F}$ are called flows, while the elements of $\mathcal{E}$ are called efforts. Those are port variables and their combination gives the power flowing inside the system. The space $\mathcal{B} = \mathcal{F} \times \mathcal{E}$ is called the bond space of power variables. Therefore the power is defined as  $\left\langle \bm{e}, \bm{f} \right\rangle = \bm{e}(\bm{f})$, where $\left\langle \bm{e} , \bm{f} \right\rangle$ is the dual product between $\bm{f}$ and $\bm{e}$.
		\begin{definition}[\cite{CourantDiracStructure}, Def. 1.1.1]
			Given the space $\mathcal{F}$ and its dual $\mathcal{E}$ with respect to the inner product $\left\langle \cdot , \cdot \right\rangle : \mathcal{F} \times \mathcal{E} \rightarrow \mathbb{R}$, define the symmetric bilinear form:
			\begin{equation}
			\left\langle \left\langle (\bm{f}_1, \bm{e}_1), (\bm{f}_2, \bm{e}_2) \right\rangle \right\rangle := \left\langle \bm{e}_1, \bm{f}_2 \right\rangle +  \left\langle \bm{e}_2, \bm{f}_1 \right\rangle, \quad \text{with} \quad (\bm{f}_i, \bm{e}_i) \in \mathcal{B}, \; i = 1, 2
			\end{equation}
			
			A Dirac structure on $\mathcal{B} := \mathcal{F} \times \mathcal{E}$ is a subspace $\mathcal{D} \subset \mathcal{B}$, {which is maximally isotropic under $\left\langle \left\langle \cdot, \cdot \right\rangle \right\rangle$.}	
		\end{definition}
		This definition can be extended to consider distributed forces and dissipation \cite{Villegas}.
		\begin{proposition}
			\label{prop:Dirac}
			Consider the space of power variables $\mathcal{F} \times \mathcal{E}$ and let $\mathcal{X}$ denote an $n$-dimensional space, the space of energy variables. Suppose that $\mathcal{F} := (\mathcal{F}_s, \ \mathcal{F}_e )$ and that $\mathcal{E} := (\mathcal{E}_s,  \ \mathcal{E}_e )$, with $\text{dim} \, \mathcal{F}_s = \text{dim} \, \mathcal{E}_s = n$ and $\text{dim} \, \mathcal{F}_e = \text{dim} \, \mathcal{E}_e = m$. Moreover, let $\bm{J}(\bm{x})$ denote
			a skew-symmetric matrix of dimension $n$ and by $\bm{B}(\bm{x})$ a matrix of dimension $n \times m$. Then, the set
			\begin{equation}
			\mathcal{D} := \left\{ (\bm{f}_s, \bm{f}_e , \bm{e}_s ,\bm{e}_e ) \in \mathcal{F} \times \mathcal{E} \vert \quad \bm{f}_s = - \bm{J}(\bm{x}) \bm{e}_s - \bm{B}(\bm{x}) \bm{f}_e, \; \bm{e}_e = \bm{B}(\bm{x})^T \bm{e}_s \right\}
			\end{equation}
			is a Dirac structure.
		\end{proposition}
		\subsection{Finite-dimensional PHs}
		Consider the time-invariant dynamical system:
		\begin{equation}
		\label{eq:finitePH}
		\begin{cases}
		\dot{ \bm{x} } &= J(\bm{x}) \nabla H(\bm{x}) + B(\bm{x})\bm{u}, \\
		\bm{y} &= B(\bm{x})^T \nabla H(\bm{x}),
		\end{cases}
		\end{equation}
		where $ H(\bm{x}) : \mathcal{X} \rightarrow \mathbb{R} $, the Hamiltonian, is a real-valued function bounded from below. Such system is called port-Hamiltonian, as it arises from the Hamiltonian modelling of a physical system and it interacts with the environment through the input $\bm{u}$ included in the formulation. The connection with the concept of Dirac structure is achieved by considering the following port behavior:
		\begin{equation}
		\begin{aligned}
		\bm{f}_s &= - \dot{ \bm{x} }, \qquad 
		&\bm{e}_s &= \diffp{H}{\bm{x}}, \\
		\bm{f}_e &= \bm{u}, \qquad
		&\bm{e}_e &= \bm{y}. \\
		\end{aligned}
		\end{equation}
		With this choice of the port variables system \eqref{eq:finitePH} defines, by proposition \ref{prop:Dirac}, a Dirac structure. Dissipation and distributed forces can be included and the corresponding system defines an extended Dirac structure, once the proper port variables are introduced.
		\subsection{Constant matrix differential operators}
		The treatment provided in \cite{MacchelliModelling} is here recovered, in order to have the most suitable definition of Stokes-Dirac structure for the mechanical systems considered in this paper.
		Let $\Omega$ denote a compact subset of $\mathbb{R}^d$ representing the spatial domain of the distributed parameter system. Then, let $\mathcal{U}$ and $\mathcal{V}$ denote two sets of smooth functions from $\Omega$ to $\mathbb{R}^{q_u}$ and $\mathbb{R}^{q_v}$ respectively.
		\begin{definition}
			A constant matrix differential operator of order $N$ is a map $L$ from $\mathcal{U}$ to $\mathcal{V}$ such that, given $\bm{u} = (u_1 , \dots , u_{q_u}) \in \mathcal{U}$ and $\bm{v} = (v_1 , . . . , v_{q_v}) \in \mathcal{V}$:
			\begin{equation}
			\label{eq:diffOp}
			\bm{v} = L \bm{u} \iff \bm{v} := \sum_{|\alpha|=0}^N  \bm{P}_{\alpha} D^{\alpha} \bm{u},
			\end{equation}
			where $\alpha := (\alpha_1, \dots , \alpha_d)$ is a multi-index of order $|\alpha| := \sum_{i=1}^d \alpha_i$, $\bm{P}_\alpha$ are a set of constant real $q_v \times q_u$ matrices and $D^{\alpha} := \partial_{x_1}^{\alpha_1} \dots \partial_{x_d}^{\alpha_d}$ is a differential operator of order $|\alpha|$ resulting from a combination of spatial derivatives. 
		\end{definition}
		\begin{definition}
			Consider the constant matrix differential operator \eqref{eq:diffOp}. Its formal adjoint is the map $L^*$ from $\mathcal{V}$ to $\mathcal{U}$ such that:
			\begin{equation}
			\bm{u} = L^* \bm{v} \iff \bm{u} := \sum_{|\alpha|=0}^N  (-1)^{|\alpha|} \bm{P}_{\alpha}^T D^{\alpha} \bm{v}.
			\end{equation}
		\end{definition}
		\begin{definition}
			\label{def:skewOp}
			Let $J$ denote a constant matrix differential operator. Then, J is skew-symmetric if and only if $J = -J^*$. This corresponds to the condition:
			\begin{equation}
			\bm{P}_{\alpha} = (-1)^{|\alpha| + 1} \bm{P}_{\alpha}^T, \quad \forall \alpha.
			\end{equation}
		\end{definition}
		An important relation between a differential operator and
		its adjoint is expressed by the following theorem.
		\begin{theorem}[\cite{PDE}, Chapter 9, theorem 9.37]
			Consider a matrix differential operator $L$ and let $L^*$ denote its formal adjoint. Then, for each function $\bm{u} \in \mathcal{U}$ and $\bm{v} \in \mathcal{V}$:
			\begin{equation}
			\int_{\Omega} \left( \bm{v}^T L \bm{u} - \bm{u}^T L^* \bm{v}\right) \d\Omega= \int_{\partial \Omega} \widetilde{B}_L(\bm{u}, \bm{v}) \d{A},
			\end{equation}
			where $\widetilde{B}_L$ is a differential operator induced on the boundary $\partial\Omega$ by L, or equivalently:
			\begin{equation}
			\bm{v}^T L \bm{u} - \bm{u}^T L^* \bm{v} = \mathrm{div} \, \widetilde{B}_L(\bm{u}, \bm{v}).
			\end{equation}
		\end{theorem}
		It is important to note that $\widetilde{B}_L$ is a constant differential operator. The quantity $\widetilde{B}_L(\bm{u}, \bm{v})$ is a constant linear combination of the functions $\bm{u}$ and $\bm{v}$ together with their spatial derivatives up to a certain order and depending on $L$.
		\begin{corollary}
			Consider a skew-symmetric differential operator $J$. Then, for each function $\bm{u} \in \mathcal{U}$ and $\bm{v} \in \mathcal{V}$ with $q_u = q_v = q$:
			\begin{equation}
			\int_{\Omega} \left( \bm{v}^T J \bm{u} + \bm{u}^T J \bm{v}\right) \d\Omega= \int_{\partial \Omega} \widetilde{B}_J(\bm{u}, \bm{v}) \d{A},
			\end{equation}
			where $\widetilde{B}_J$ is a symmetric differential operator on $\partial\Omega$ depending on the differential operator~$J$.
		\end{corollary}
		\subsection{Constant Stokes-Dirac structures}
		Following again \cite{MacchelliModelling} let $\mathcal{F}$ denote the space of flows, i.e. the space of smooth functions from the compact set $\Omega \subset \mathbb{R}^d$ to $\mathbb{R}^q$. For simplicity assume that  the space of efforts is $\mathcal{E} \equiv \mathcal{F}$ (generally speaking these spaces are Hilbert spaces linked by duality, as in \cite{Villegas}). Given $\bm{f} = (f_1, \dots, f_q) \in \mathcal{F}$ and $\bm{e} = (e_1, \dots, e_q) \in \mathcal{E}$. Let $\bm{z} = B_\partial(\bm{e})$ denote the boundary terms, where $B_\partial$ provides the restriction on $\partial\Omega$ of the effort $e$ and of its spatial derivatives of proper order. Then it can be written:
		\begin{equation}
		\int_{\partial \Omega} \widetilde{B}_J(\bm{e}_1, \bm{e}_2) \d{A} = \int_{\partial \Omega} B_J(\bm{z}_1, \bm{z}_2) \d{A}, \quad \text{with} \quad  \widetilde{B}_J(\cdot, \cdot) = B_J(B_{\partial}(\cdot), \, B_{\partial}(\cdot)).
		\end{equation}
		Define the set
		\begin{equation}
		\mathcal{Z} := \left\{ \bm{z} \vert \bm{z} = B_{\partial}(\bm{e})  \right\}.
		\end{equation}
		\begin{theorem}[\cite{MacchelliModelling}]
			\label{th:StokesDirac}
			Consider the space of power variables $\mathcal{B} = \mathcal{F} \times \mathcal{E} \times \mathcal{Z}$. The linear subspace $\mathcal{D} \subset \mathcal{B}$
			\begin{equation}
			\mathcal{D} = \left\{ (\bm{f}, \bm{e}, \bm{z}) \in  \mathcal{F} \times \mathcal{E} \times \mathcal{Z} \; \vert \; \bm{f} = -J \bm{e}, \; z = B_\partial(\bm{e}) \right\},
			\end{equation}
			is a Stokes–Dirac structure on $\mathcal{B}$ with respect to the pairing
			\begin{equation}
			\left\langle \left\langle (\bm{f}_1, \bm{e}_1, \bm{z}_1), (\bm{f}_2, \bm{e}_2, \bm{z}_2) \right\rangle \right\rangle  := \int_{\Omega} \left( \bm{e}_1^T \bm{f}_2 + \bm{e}_2^T \bm{f}_1 \right) \d\Omega + \int_{\partial \Omega} B_J(\bm{z}_1, \bm{z}_2) \d{A}.
			\end{equation}
		\end{theorem}
		\begin{remark}
			The constant Stokes-Dirac structure has been defined in case of smooth vector valued functions for simplicity. In this context the pairing has been defined as the $L^2$ inner product of vector-valued function. The definition is indeed more general and encompasses the case of more complex functional spaces. The Mindlin plate for example is defined on a mixed function space of scalar-, vector- and tensor- valued functions. The result presented here remains valid provided that the proper {pairing is being chosen}. Furthermore, the constant differential operator may contain intrinsic operators ($\mathrm{div}, \, \mathrm{grad}$) as it will be shown for the Mindlin plate case.
		\end{remark}
	}
	
	\subsection{Infinite-dimensional PHs}
	Following \cite[Chapter~3]{BookZwart}, the prototypical example of the Timoshenko beam will be used to illustrate the class of distributed port-Hamiltonian systems. This model consists of two coupled PDEs, describing the vertical displacement and rotation scalar fields:
	\begin{equation}
	\label{eq:TimoModel}
	\begin{cases}
	\displaystyle \rho(x) \diffp[2]{w}{t}(x,t) &= \displaystyle \diffp{}{x} \left[K(x) \left(\diffp{w}{x}(x,t) - \phi(x,t)\right)\right], \quad x \in (0,L),\, t \ge 0 \\
	\displaystyle I_{\rho}(x) \diffp[2]{\phi}{t}(x,t) &= \displaystyle \diffp{}{x} \left(EI(x) \diffp{ \phi}{x}(x,t)\right) + K(x) \left(\diffp{w}{x}(x,t) - \phi(x,t) \right),
	\end{cases}
	\end{equation}
	where ${w}(x,t)$ is the transverse displacement and $\phi(x,t)$ is the rotation angle of a fiber of the beam. The coefficients $\rho(x), I_{\rho}(x), E(x), I(x)$ and $K(x)$ are the mass per unit length, the rotary moment of inertia of a cross section, Young's modulus of elasticity, the moment of inertia of a cross section and the shear modulus, respectively. The energy variables are chosen as follows:
	\begin{equation}
	\begin{aligned}
	\alpha_{w} &:= \rho(x) \diffp{w}{t}(x,t), \quad &\text{{Linear Momentum,}} \\
	\alpha_{\phi} &:= I_{\rho}(x) \diffp{\phi}{t}(x,t), \quad &\text{{Angular Momentum,}} \\
	\alpha_{\kappa} &:= \diffp{\phi}{x}(x,t), \quad &\text{{Curvature,}} \\
	\alpha_{\gamma} &:= \diffp{w}{x}(x,t) - \phi(x,t), \quad &\text{{Shear Deformation}}. \\
	\end{aligned}
	\end{equation}
	
	Those variables are collected in the vector $\bm{\alpha} := (\alpha_{w}, \, \alpha_{\phi}, \, \alpha_{\kappa}, \, \alpha_{\gamma} )^T $, so that the Hamiltonian can be written as a quadratic functional in the energy variables: 
	\begin{equation}
	H = \frac{1}{2} \int_{0}^{L} \bm{\alpha}^T Q \bm{\alpha} \; dx,
	\qquad \text{where} \qquad
	Q = 
	\begin{bmatrix}
	\frac{1}{\rho(x)} & 0 & 0 & 0 \\
	0 & \frac{1}{I_{\rho}(x)} & 0 & 0 \\
	0 & 0 & EI(x) & 0 \\
	0 & 0 & 0 & K(x) \\
	\end{bmatrix}.
	\end{equation}
	
	The co-energy variables are found by computing the variational derivative of the Hamiltonian (see  \cite{MacchelliTimo}):
	\begin{equation}
	\begin{aligned}
	e_{w} &:= \diffd{H}{\alpha_w} = \diffp{w}{t}(x,t), \quad &\text{{Linear Velocity,}} \\
	e_{\phi} &:= \diffd{H}{\alpha_{\phi}} = \diffp{\phi}{t}(x,t), \quad &\text{{Angular Velocity,}} \\
	e_{\kappa} &:= \diffd{H}{\alpha_{\kappa}} =EI(x) \diffp{\phi}{x}(x,t), \quad &\text{{Flexural Momentum,}}\\
	e_{\gamma} &:= \diffd{H}{\alpha_{\gamma}} = K(x) \left(\diffp{w}{x}(x,t) - \phi(x,t) \right), \quad &\text{{Shear Force}}. \\
	\end{aligned}
	\end{equation}
	These variables are again collected in the vector $\bm{e} = (e_{w}, \, e_{\phi}, \, e_{\kappa}, \, e_{\gamma} )^T $, so that {system \eqref{eq:TimoModel} can be rewritten in terms of energy and co-energy variables:}
	\begin{equation}
	\label{eq:PH_Timo}
	\diffp{\bm{\alpha}}{t} = J \mathbf{e},  	\qquad \text{where} \qquad
	J = 
	\begin{bmatrix}
	0 & 0 & 0 & \diffp{}{x} \\
	0 & 0 & \diffp{}{x} & 1  \\
	0 & \diffp{}{x} & 0 & 0 \\
	\diffp{}{x} & -1 & 0 & 0 \\
	\end{bmatrix} .
	\end{equation}
	{By Definition \eqref{def:skewOp}, the constant matrix differential operator $J$ is skew-symmetric. In order to unveil the Stokes-Dirac structure for this model, boundary variables have to be defined. The energy rate is given by (see \cite{BookZwart} for a detailed derivation):
		\begin{equation}
		\dot{H} = \bm{f}_{\partial}^T \bm{e}_{\partial},
		\end{equation}
		where 
		\begin{equation}
		\begin{aligned}
		\bm{f}_{\partial} &= 
		\begin{bmatrix}
		e_w(0) & e_{\phi}(0) & e_{\gamma}(L) & e_{\kappa}(L) \\
		\end{bmatrix}^T, \\
		\bm{e}_{\partial} &= 
		\begin{bmatrix}
		-e_{\gamma}(0) & -e_{\kappa}(0) & e_{w}(L) & e_{\phi}(L) \\
		\end{bmatrix}^T.
		\end{aligned}
		\end{equation}
		Consider now the power space
		\begin{equation}
		\label{eq:bondTimo}
		\mathcal{B} := \left\{(\bm{f}, \bm{e}, \bm{z}) \in \mathcal{F} \times \mathcal{E} \times \mathcal{Z} \right\},
		\end{equation}
		where $\mathcal{F} \equiv \mathcal{E} =  L^2( (0, L); \mathbb{R}^4)$ and
		\begin{equation}
		\mathcal{Z} := \left\{ \bm{z} \; \vert \; \bm{z} = \begin{pmatrix} \bm{f}_{\partial} \\ \bm{e}_{\partial} \end{pmatrix} \right\}.
		\end{equation}
		The duality pairing between elements of $\mathcal{B}$ is then defined as follows:
		\begin{equation}
		\label{eq:bil_Timo}
		\left\langle \left\langle (\bm{f}_1, \bm{e}_1, \bm{z}_1), (\bm{f}_2, \bm{e}_2, \bm{z}_2) \right\rangle \right\rangle :=  \int_{0}^L \left( \bm{e}_1^T \bm{f}_2 + \bm{e}_2^T \bm{f}_1 \right) \d{x} + B_J(\bm{z}_1, \bm{z}_2) , 
		\end{equation}
		where $B_J(\bm{z}_1, \bm{z}_2) := (\bm{f}_{\partial, 1})^T \bm{e}_{\partial, 2} + (\bm{f}_{\partial, 2}) ^T \bm{e}_{\partial, 1}$. It is now possible to state the following theorem:
		\begin{theorem}[Stokes-Dirac Structure for the Timoshenko beam]
			Consider the space of power variables $\mathcal{B}$ defined in \eqref{eq:bondTimo}. By theorem \ref{th:StokesDirac} the linear subspace $\mathcal{D} \subset \mathcal{B}$
			\begin{equation}
			\mathcal{D} =  \left\{(\bm{f}, \bm{e},\bm{z}) \in \mathcal{B} | \; \bm{f} = - \diffp{\bm{\alpha}}{t} = -J\bm{e}, \; \bm{z} = \begin{pmatrix} \bm{f}_{\partial} \\ \bm{e}_{\partial} \end{pmatrix} 
			\right\},
			\end{equation}
			is a Stokes-Dirac structure with respect to the pairing $\left\langle \left\langle \cdot, \cdot \right\rangle \right\rangle$ given by \eqref{eq:bil_Timo}.
		\end{theorem}
	}
	
	\section{Mindlin theory for thick plates}
	\label{sec:Min_Var}
	{In this section the classical model of thick plates is first recalled by making use of a tensorial formalism to simplify the discussion on the port-Hamiltonian formulation.} The Mindlin plate theory, originally presented in \cite{mindlin}, is more suited for plates having a large thickness. The fibers, i.e. a segment perpendicular to the mid-plane, of the plate are supposed to remain straight after the deformation, but not necessarily normal to the mid-plane. For this reason two new kinematic variables have to be added, in order to represent the deflection of the cross sections. { Let $\theta_x$ denote the deflection of the cross section with respect to the opposite side of the $y$ axis and  $\theta_y$ its deflection along the $x$ axis. The displacement field is therefore given by the following relations:
		\begin{equation}
		\begin{aligned}
		u(x,y,z) &= -z \theta_x(x,y), \quad &\text{Displacement along $x$}, \\
		v(x,y,z) &= -z \theta_y(x,y), \quad &\text{Displacement along $y$}, \\
		w(x,y,z) &= w(x,y),  \quad &\text{Displacement along $z$}. \\
		\end{aligned}
		\end{equation}
		Both bending and shear deformations are taken into account. The bending part is described by the second order shear tensor:
		\begin{equation}
		\mathbb{E}_b = 
		\begin{bmatrix}
		\epsilon_{xx} & \gamma_{xy} \\
		\gamma_{xy} & \epsilon_{yy} \\
		\end{bmatrix} = -z \; \mathrm{Grad}(\bm\theta) = -z
		\begin{bmatrix}
		\diffp{\theta_x}{x} & \frac{1}{2} \left( \diffp{\theta_y}{x} + \diffp{\theta_x}{y} \right)\\
		\frac{1}{2} \left( \diffp{\theta_y}{x} + \diffp{\theta_x}{y} \right) & \diffp{\theta_y}{y} \\
		\end{bmatrix},
		\end{equation}
		where $\bm{\theta} = (\theta_x, \, \theta_y)^T$. The tensor $\mathrm{Grad}(\bm{\theta})$ stands for the symmetric gradient of the vector $\bm{\theta}$ and corresponds to the curvature tensor:
		\begin{equation}
		\mathbb{K} = \begin{bmatrix}
		\kappa_{xx} &  \kappa_{xy}\\
		\kappa_{xy} & \kappa_{yy} \\
		\end{bmatrix} := \mathrm{Grad}(\bm{\theta}) =  \frac{1}{2} \left(\nabla \otimes \bm{\theta} + \left(\nabla \otimes \bm{\theta}\right)^T \right),
		\end{equation}
		where $\bm{u} \otimes {\bm{v}}$ denotes the outer product of vectors, equivalent to the dyadic product given by $\bm{u}\bm{v}^T$. The corresponding bending  stress field $\mathbb{S}_b$ is obtained by considering the constitutive relation which, for an isotropic homogeneous material, reads:
		\begin{equation}
		\mathbb{S}_{b} = \frac{E}{1-\nu^2} \left\{ (1-\nu) \mathbb{E}_b +\nu \bm{I}_{2 \times 2} \Tr(\mathbb{E}_b) \right\},
		\end{equation}
		where $E$ is Young's modulus, $\nu$ Poisson's ratio, $\bm{I}_{2 \times 2}$ the identity operator in $\mathbb{R}^2$ and $\Tr$ the trace operator. By averaging the torques produced by stresses along a fiber, it is possible to define the flexural momenta tensor:
		\begin{equation}
		\mathbb{M}_{ij} = \int_{-\frac{h}{2} }^{ \frac{h}{2} } -z \, \mathbb{S}_{b, \; ij} \d z = \mathbb{D}_{ijkl} \,  \mathbb{K}_{kl}, 
		\end{equation}
		where $h$ is the plate thickness and $\mathbb{D}$ is the a fourth order symmetric tensor and represents the bending rigidity tensor. The components of the momenta tensor are given by the following relations:
		\begin{equation}
		\mathbb{M} = \begin{bmatrix}
		M_{xx} &  M_{xy}\\
		M_{xy} &  M_{yy} \\
		\end{bmatrix},  \qquad 
		\begin{aligned}
		M_{xx} &= D\left(\kappa_{xx} + \nu \kappa_{yy} \right),\\
		M_{yy} &= D\left(\kappa_{yy} + \nu \kappa_{xx} \right),\\
		M_{xy} &= D\left(1 - \nu \right) \kappa_{xy}, \\
		\end{aligned} 
		\end{equation}
		where $D = \frac{E h^3}{12 (1 - \nu^2)}$ is the bending rigidity. The shear deformations are described using the shear strain vector:
		\begin{equation}	
		\bm{\epsilon}_s = 
		\begin{pmatrix}
		\gamma_{xz} \\
		\gamma_{yz} \\
		\end{pmatrix} = 
		\begin{pmatrix}
		\diffp{w}{x} - \theta_x \\
		\diffp{w}{y} - \theta_y \\
		\end{pmatrix}.
		\end{equation}
		The corresponding stress field is simply given by $\bm\sigma_s = G \bm{\epsilon}_s$, where $G$ is the shear modulus $G = \frac{E}{2 (1 + \nu)}$. The shear forces are again obtained by averaging the shear stress along plate fibers. However, given the excessive rigidity of the shear contribution, a correction factor $k = 5/6$ (see \cite{mindlin}) is introduced:
		\begin{equation}
		\bm{q} = \begin{pmatrix}
		q_x \\
		q_y \\
		\end{pmatrix} =
		\int_{-\frac{h}{2} }^{ \frac{h}{2} } k \bm\sigma_s \d{z} = k G h \, \bm{\epsilon}_s = k G h
		\begin{pmatrix}
		\gamma_{xz} \\
		\gamma_{yz} \\
		\end{pmatrix}.
		\end{equation}
		All the needed variables are now defined and the equations of motions can be recalled (see \cite{mindlin}): 
		\begin{equation}
		\begin{cases}
		\displaystyle \rho h \diffp[2]{w}{t} &= \mathrm{div}(\bm{q}),  \vspace{1mm}\\
		\displaystyle \rho\frac{h^3}{12} \diffp[2]{\bm \theta}{t} &= \bm{q} + \mathrm{Div}(\mathbb{M}), \\
		\end{cases}
		\end{equation}
		where $\rho$ is the plate mass density. The differential operators $\mathrm{div}()$ and $\mathrm{Div}()$ denote the divergence of a vector and of a tensor, respectively:
		\begin{equation*}
		\bm{a} = \mathrm{Div}(\mathbb{A})  \qquad \text{with } \bm{a}_i = \mathrm{div}(\mathbb{A}_{ji}) = \sum_{j = 1}^n \diffp{\mathbb{A}_{ji}}{x_j},
		\end{equation*}
		where $\bm{a}$ is defined column-wise. The kinetic and potential energy densities per unit area ($\mathcal{K}$ and $\mathcal{U}$), are given respectively by:
		\begin{align*}
		\mathcal{K} &= \frac{1}{2} \rho \left\{ h \left(\diffp{w}{t}\right)^2 +  \frac{h^3}{12} \diffp{\bm\theta}{t} \cdot \diffp{\bm\theta}{t} \right\} ,\\
		\mathcal{U} &= \frac{1}{2} \left\{ \mathbb{M} \cddot \mathbb{K} + \bm{q} \cdot \bm{\epsilon}_s \right\} .	
		\end{align*} 
		where the tensor contraction $ \mathbb{M} \cddot \mathbb{K}$ in Cartesian coordinates is expressed as:
		\[\mathbb{M} \cddot \mathbb{K} = \sum_{i,j = 1}^{2} M_{ij} \kappa_{ij} = \Tr(\mathbb{M}^T \mathbb{K}). \]
		The total energy density is the sum of kinetic and potential energies
		\begin{equation}
		\mathcal{H} = \mathcal{K} + \mathcal{U}
		\end{equation}
		and the corresponding energies
		\begin{equation}
		H = \int_{\Omega} \mathcal{H} \ \d\Omega, \qquad K = \int_{\Omega} \mathcal{K} \ \d\Omega, \qquad U = \int_{\Omega} \mathcal{U} \ \d\Omega,
		\end{equation}
		where $\Omega$ is an open connected set of $\mathbb{R}^2$.
	}
	
	\section{PH tensorial formulation of the Mindlin plate}
	\label{sec:PH_ten_Min}
	In this section the new tensorial formulation for the Mindlin plate is presented. {Obviously this result is equivalent to the one found in the vectorial case \cite{MacchelliMindlin}, but the tensorial formulation is more suitable from the point of view of functional analysis since it makes appear intrinsic operators ($\mathrm{div}, \mathrm{Div}, \mathrm{grad}, \mathrm{Grad}$) as it will be shown in the following, regardless of the choice of coordinate frame. \newline
		The Hamiltonian energy is first considered in order to find the proper energy variables:
		\begin{equation}
		\label{eq:H_min}
		H = \int_{\Omega} \frac{1}{2} \left\{ \rho h \left(\diffp{w}{t} \right)^2 + \frac{\rho h^3}{12} \diffp{\bm{\theta}}{t} \cdot   \diffp{\bm{\theta}}{t} +   \mathbb{M} \cddot \mathbb{K} + \bm{q} \cdot \bm{\epsilon}_s  \right\}  \d\Omega, 
		\end{equation}
		The choice of the energy variables is the same as in \cite{MacchelliMindlin} but here scalar,}  vector and tensor variables are gathered together:
	\begin{equation}
	\begin{aligned}
	\alpha_w &= \rho h \diffp{w}{t}, \quad &\text{Linear momentum,} \\
	\mathbb{A}_{\kappa} &= \mathbb{K}, \quad &\text{Curvature Tensor,} \\
	\end{aligned} \qquad
	\begin{aligned}
	\bm\alpha_{\theta} &=  \frac{\rho h^3}{12} \diffp{\bm{\theta}}{t}, \quad &\text{Angular momentum,}\\
	\bm\alpha_{\gamma} &= \bm{\epsilon}_s. \quad &\text{Shear Deformation.}\\
	\end{aligned}
	\end{equation}
	The co-energy variables are found by computing the variational derivative of the Hamiltonian:
	\begin{equation}
	\begin{aligned}
	e_w &:= \diffd{H}{\alpha_w} = \diffp{w}{t},  \quad &\text{Linear velocity,} \\
	\mathbb{E}_{\kappa} &:= \diffd{H}{\mathbb{A}_{\kappa}} = \mathbb{M}, \quad &\text{Momenta Tensor,}\\
	\end{aligned} \qquad
	\begin{aligned}
	\bm{e}_{\theta} &:= \diffd{H}{\bm\alpha_{\theta}} = \diffp{\bm{\theta}}{t}, \quad &\text{Angular velocity,}  \\
	\bm{e}_{\gamma} &:= \diffd{H}{\bm{\epsilon}_s} = \bm{q} \quad &\text{Shear stress.} \\
	\end{aligned}
	\end{equation}
	\begin{proposition}
		The variational derivative of the Hamiltonian with respect to the curvature tensor is the momenta tensor $\diffd{H}{\mathbb{A}_{\kappa}} = \mathbb{M}$.
		\begin{proof}
			The contribution due to the bending part in Hamiltonian is given by:
			\[H_b(\mathbb{K}) = \frac{1}{2} \int_{\Omega}  \mathbb{M} \cddot \mathbb{K} \; \d\Omega = \frac{1}{2} \int_{\Omega} \Tr(\mathbb{M}^T \mathbb{K}) \; \d\Omega, \]
			where the momenta tensor depends on the curvatures tensor $\mathbb{M}_{ij} = \mathbb{D}_{ijkl}\mathbb{K}_{kl} $ where $\mathbb{D} = \mathbb{D}^T$ is a fourth order symmetric positive definite tensor.
			So a variation $\delta\mathbb{K}$ of the curvatures tensor with respect to a given value $\mathbb{K}_0$ leads to:
			\[H_b(\mathbb{K}_0+ \varepsilon \delta\mathbb{K}) = \frac{1}{2} \int_{\Omega} \Tr(\mathbb{M}_0^T \mathbb{K}_0) \; \d\Omega + \varepsilon \frac{1}{2} \int_{\Omega} \left\{ \Tr(\mathbb{M}_0^T \delta\mathbb{K}) + \Tr(\delta\mathbb{M}^T \mathbb{K}_0)\right\}\; \d\Omega  + O(\varepsilon^2). \]
			The term $\Tr(\delta\mathbb{M}^T \mathbb{K}_0)$ can be further manipulated as follows 
			\[ \Tr(\delta\mathbb{M}^T \mathbb{K}_0) = \Tr(\mathbb{K}_0^T \delta\mathbb{M}) = \Tr(\mathbb{K}_0^T \, \mathbb{D} \, \delta\mathbb{K}) = \Tr(\mathbb{M}_0^T \delta\mathbb{K}),
			\]
			so that 
			\[H_b(\mathbb{K}_0+ \varepsilon \delta\mathbb{K}) = \frac{1}{2} \int_{\Omega} \Tr(\mathbb{M}_0^T \mathbb{K}_0) \; \d\Omega + \varepsilon \int_{\Omega} \left\{ \Tr(\mathbb{M}_0^T \delta\mathbb{K})\right\}\; \d\Omega  + O(\varepsilon^2). \]
			By definition of the variational derivative (see e.g. \cite{VANDERSCHAFT2002166}) it can be written:
			\[H_b(\mathbb{K}_0+ \varepsilon \delta\mathbb{K}) = H_b(\mathbb{K}_0) + \varepsilon \left\langle \diffd{H_b}{\mathbb{K}} , \delta\mathbb{K}  \right\rangle_{\mathscr{H}} + O(\varepsilon^2), \]
			where $\mathscr{H}$ is the space of the square integrable symmetric tensors endowed with the integral of the tensor contraction as inner product. Then, by identification $\diffd{H_b}{\mathbb{K}} = \diffd{H}{\mathbb{K}} = \mathbb{M}_0$.
		\end{proof}
	\end{proposition}
	
	The port-Hamiltonian system is expressed as follows:
	\begin{equation}
	\begin{cases}
	\displaystyle\diffp{\alpha_w}{t} &= \mathrm{div}(\bm{e}_{\gamma}), \vspace{1mm} \\
	\displaystyle\diffp{\bm\alpha_\theta}{t} &= \mathrm{Div}( \mathbb{E}_{\kappa}) + \bm{e}_{\gamma}, \vspace{1mm} \\
	\displaystyle\diffp{\mathbb{A}_\kappa}{t} &= \mathrm{Grad}(\bm{e}_{\theta}), \vspace{1mm} \\
	\displaystyle\diffp{ \bm\alpha_{\gamma} }{t} &= \mathrm{grad}(e_w) - \bm{e}_{\theta},
	\end{cases}
	\end{equation}
	
	If the variables are concatenated together, the formally skew-symmetric operator $J$ can be highlighted:
	\begin{equation}
	\label{eq:PH_sys_Min_Ten}
	\diffp{}{t}
	\begin{pmatrix}
	\alpha_w \\
	\bm\alpha_\theta \\
	\mathbb{A}_\kappa \\
	\bm\alpha_{\gamma} \\
	\end{pmatrix} = 
	\underbrace{\begin{bmatrix}
		0  & 0  & 0  & \mathrm{div} \\
		0 & 0 &  \mathrm{Div} & \bm{I}_{2 \times 2}\\
		0  & \mathrm{Grad}  & 0  & 0\\
		\mathrm{grad} & -\bm{I}_{2 \times 2} &  0 & 0  \\
		\end{bmatrix}}_{J}
	\begin{pmatrix}
	e_w \\
	\bm{e}_{\theta} \\
	\mathbb{E}_{\kappa} \\
	\bm{e}_{\gamma} \\
	\end{pmatrix},
	\end{equation}
	where all zeros are intended as nullifying operators from the space of input variables to the space of output variables.
	\begin{remark}
		It can be observed that the interconnection structure given by $J$ in \eqref{eq:PH_sys_Min_Ten} mimics that of the Timoshenko beam given in \eqref{eq:PH_Timo}.
	\end{remark}

	\begin{theorem}{The adjoint of the tensor divergence $\mathrm{Div}$ is $- \mathrm{Grad}$, the opposite of the symmetric gradient.}
		\begin{proof}
			Let us consider the Hilbert space of the square integrable symmetric tensors of size $n \times n$ over an open connected set $\Omega$. This space will be denoted by $\mathscr{H}_1 = L^2(\Omega, \mathbb{R}^{n \times n}_{\text{sym}})$. This space is endowed with the integral of the tensor contraction as scalar product:
			\[\left\langle \mathbb{E} , \mathbb{F} \right\rangle_{\mathscr{H}_1} = \int_{\Omega}  \mathbb{E} \cddot \mathbb{F} \; \d\Omega = \int_{\Omega} \Tr(\mathbb{E}^T \mathbb{F}) \; \d\Omega, \quad \forall \, \mathbb{E} , \mathbb{F} \in L^2(\Omega, \mathbb{R}^{n \times n}_{\text{sym}}). \]
			
			Moreover the Hilbert space of the square integrable vector functions over the same open connected set $\Omega$ will be denoted by $\mathscr{H}_2 = L^2(\Omega, \mathbb{R}^{n})$. This space is endowed with the following scalar product:
			\[\left\langle \bm{\varepsilon} , \bm{\phi} \right\rangle_{\mathscr{H}_2} = \int_{\Omega}  \bm{\varepsilon} \cdot \bm{\phi} \; \d\Omega = \int_{\Omega} \bm{\varepsilon}^T \bm{\phi}\; \d\Omega, \quad \forall \bm{\varepsilon}, \bm{\phi} \in L^2(\Omega, \mathbb{R}^{n}). \]
			Let us consider the tensor divergence operator defined as:
			\[
			\begin{aligned}
			A: \; \mathscr{H}_1& \rightarrow \mathscr{H}_2, \\
			\mathbb{E}& \rightarrow \mathrm{Div}(\mathbb{E}) = \bm{\varepsilon}, \\
			\end{aligned}
			\qquad \text{with } \bm{\varepsilon}_i = \mathrm{div}(\mathbb{E}_{ji}) = \sum_{j = 1}^n \diffp{\mathbb{E}_{ji}}{x_j}.
			\]
			We try to identify $A^*$
			\[
			\begin{aligned}
			A^*: \; \mathscr{H}_2& \rightarrow \mathscr{H}_1, \\
			\bm{\phi}& \rightarrow  A^* \bm{\phi} = \mathbb{F}, \\
			\end{aligned}
			\]
			such that \[
			\left\langle A \mathbb{E} , \bm{\phi} \right\rangle_{\mathscr{H}_2} = \left\langle \mathbb{E} , A^* \bm{\phi} \right\rangle_{\mathscr{H}_1},
			\begin{aligned} \qquad
			&\forall \mathbb{E} \in \mathrm{Domain}(A) \subset \mathscr{H}_1 \\
			&\forall \bm{\phi} \in \mathrm{Domain}(A^*) \subset \mathscr{H}_2
			\end{aligned}
			\]
			Now let us take $\mathbb{E} \in \mathcal{C}_0^1(\Omega, \mathbb{R}^{n \times n}_{\text{sym}}) \subset \mathrm{Domain}(A)$ the space of differentiable symmetric tensors with compact support in~$\Omega$. Additionally $\bm{\phi}$ will belong to $\mathcal{C}_0^1(\Omega, \mathbb{R}^n) \subset \mathrm{Domain}(A^*)$, the space of differentiable vector functions with compact support in~$\Omega$. Then
			\[
			\begin{aligned}
			\left\langle \mathrm{Div}\mathbb{E} , \bm{\phi} \right\rangle_{\mathscr{H}_2} &= \int_{\Omega}  \bm{\varepsilon} \cdot \bm{\phi} \; \d\Omega\, \\
			&= \int_{\Omega} \sum_{i=1}^n \sum_{j=1}^n \diffp{\mathbb{E}_{ji}}{x_j} \phi_i \; \d\Omega\,,  \\ 
			&= - \int_{\Omega} \sum_{i=1}^n \sum_{j=1}^n \mathbb{E}_{ji} \diffp{\phi_i}{x_j} \; \d\Omega\,, \qquad &\text{since the functions vanish at the boundary,}\\
			&= - \int_{\Omega} \sum_{i=1}^n \sum_{j=1}^n \mathbb{E}_{ji} \mathbb{F}_{ji} \;\d\Omega\,,  \qquad &\text{if we first choose} \quad \mathbb{F}_{ji} = \diffp{\phi_i}{x_j}\,.\\
			\end{aligned}	 
			\]
			But in this latter case, it could not indeed  be stated that $\mathbb{F} \in L^2(\Omega, \mathbb{R}^{n \times n}_{\text{sym}})$. Now, since  $\mathbb{E} \in L^2(\Omega, \mathbb{R}^{n \times n}_{\text{sym}})$, $\mathbb{E}_{ji}=\mathbb{E}_{ij}$,  thus we are  allowed to further decompose the last equality as:
			
			\[ \sum_{i,j}\mathbb{E}_{ji} \diffp{\phi_i}{x_j} = \sum_{i,j}\mathbb{E}_{ji} \frac{1}{2} \left(\diffp{\phi_i}{x_j} + \diffp{\phi_j}{x_i}  \right) = 	\sum_{i,j} \mathbb{E}_{ji} \mathbb{F}_{ji}, \qquad \text{with } \mathbb{F}_{ji}:= \frac{1}{2} \left(\diffp{\phi_i}{x_j} + \diffp{\phi_j}{x_i}  \right).
			\]
			Thus $\mathbb{F} \in L^2(\Omega, \mathbb{R}^{n \times n}_{\text{sym}})$ and it can be stated that:
			\[ \left\langle \mathrm{Div}\mathbb{E} , \bm{\phi} \right\rangle_{\mathscr{H}_2} = - \int_{\Omega} \sum_{i,j} \mathbb{E}_{ji} \frac{1}{2} \left(\diffp{\phi_i}{x_j} + \diffp{\phi_j}{x_i}  \right) \;\d\Omega = - \int_{\Omega} \sum_{i,j} \mathbb{E}_{ji} \mathbb{F}_{ji} \;\d\Omega = \left\langle \mathbb{E} , -\mathrm{Grad} \bm{\phi} \right\rangle_{\mathscr{H}_1}. \]
			It can be concluded that the formal adjoint of $\mathrm{Div}$ is $\mathrm{Div}^* = -\mathrm{Grad}$.
		\end{proof}
	\end{theorem}
	\begin{figure}[t]
		\centering
		\includegraphics[width=0.7\textwidth]{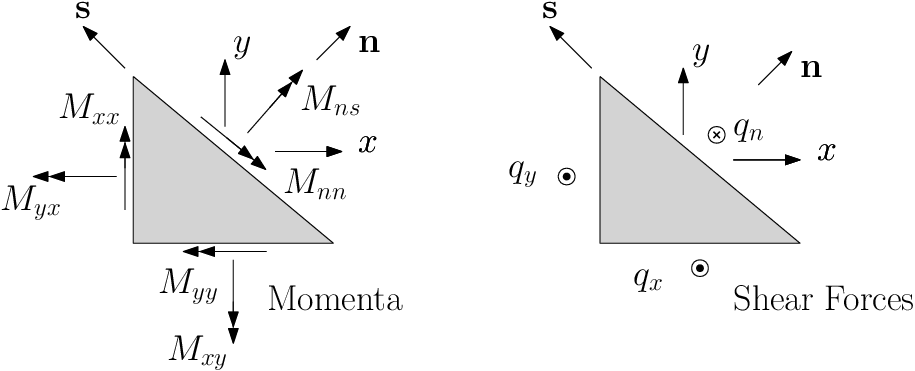}
		\caption{Cauchy law for momenta and forces at the boundary.}
		\label{fig:Cauchy_law}
	\end{figure}
	\begin{figure}[t]
		\centering
		\includegraphics[width=0.7\textwidth]{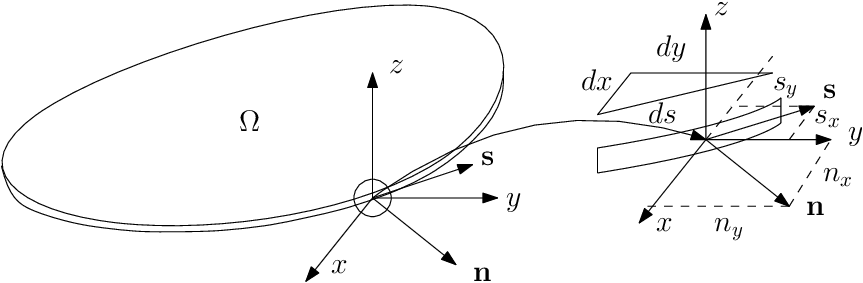}
		\caption{Reference frames and notations.}
		\label{fig:plate_ref}
	\end{figure}
	{
		We shall now establish the total energy balance in terms of boundary variables. Those will then be part of the underlying Stokes-Dirac structure of this model:
		\begin{equation}
		\label{eq:PowRate}
		\begin{aligned}
		\dot{H}&= \int_{\Omega} \left\{ \diffp{\alpha_w}{t} e_w  + \diffp{\bm\alpha_\theta}{t} \cdot \bm{e}_\theta + \diffp{\mathbb{A}_{\kappa}}{t} \cddot \mathbb{E}_{\kappa}  + \diffp{\bm\alpha_{\gamma}}{t} \cdot \bm{e}_{\gamma} \right\} \d\Omega\\
		&= \int_{\Omega} \left\{ \mathrm{div}(\bm{e}_{\gamma}) e_w  + \mathrm{Div}(\mathbb{E}_{\kappa}) \cdot \bm{e}_\theta + \; \mathrm{Grad}(\bm{e}_{\theta}) \cddot \mathbb{E}_{\kappa}  + \mathrm{grad} (e_w) \cdot \bm{e}_{\gamma} \right\} \d\Omega \\
		&= \int_{\partial \Omega} \left\{ w_t \, q_n  + \omega_n \, M_{nn} + \omega_s \, M_{ns} \right\} \d{s},  \\
		\end{aligned}
		\end{equation}
		where $s$ is the curvilinear abscissa. The last integral is obtained by applying Green-Gauss theorem. The boundary variables appearing in the last line of \eqref{eq:PowRate} and illustrated in Fig.~\ref{fig:Cauchy_law} are defined as follows:
		\begin{equation}
		\label{eq:QnMnnMns}
		\begin{aligned}
		\text{Shear Force}  \qquad q_{n} &:= \bm{q} \cdot \bm{n}=  \bm{e}_{\gamma} \cdot \bm{n},  \\
		\text{Flexural momentum} \quad 
		M_{nn} &:=  \mathbb{M} \cddot (\bm{n}\otimes{\bm{n}}) = \mathbb{E}_{\kappa} \cddot (\bm{n}\otimes{\bm{n}}) 	\\
		\text{Torsional momentum} \quad M_{ns} &:= \mathbb{M} \cddot (\bm{s}\otimes{\bm{n}}) = \mathbb{E}_{\kappa} \cddot (\bm{s}\otimes{\bm{n}}),	
		\end{aligned},
		\end{equation}
		Vectors $\bm{n}$ and $\bm{s}$ designate the normal and tangential unit vectors to the boundary, as shown in Fig. \ref{fig:plate_ref}. The corresponding power conjugated variables are:
		\begin{equation}
		\label{eq:wtwnws}
		\begin{aligned}
		\text{Vertical velocity}  \quad w_t &:= \diffp{w}{t} = e_w, \\
		\text{Flexural rotation} \quad 
		\omega_{n} &:= \diffp{\bm{\theta}}{t} \cdot \bm{n} = \bm{e}_\theta \cdot \bm{n} \\
		\text{Torsional rotation} \quad 
		\omega_{s} &:= \diffp{\bm{\theta}}{t} \cdot \bm{s} = \bm{e}_\theta \cdot \bm{s}.	
		\end{aligned},
		\end{equation}
		Analogously to the Timoshenko beam case, the Stokes-Dirac structure for the Mindlin plate can now be defined. Consider now the bond space:
		\begin{equation}
		\label{eq:bondMin}
		\mathcal{B} := \left\{(\bm{f}, \bm{e}, \bm{z}) \in \mathcal{F} \times \mathcal{E} \times \mathcal{Z} \right\},
		\end{equation}
		where $\mathcal{F}=  \mathscr{L}^2(\Omega) := 
		L^2(\Omega) \times L^2(\Omega; \mathbb{R}^2) \times L^2(\Omega, \mathbb{R}^{2 \times 2}_{\text{sym}}) \times L^2(\Omega; \mathbb{R}^2)$ and $ \mathcal{E} =  \mathscr{H}^1(\Omega) = H^{1}(\Omega) \times H^{1}(\Omega, \mathbb{R}^2) \times H^{\text{Div}}(\Omega, \mathbb{R}^{2 \times 2}_{\text{sym}}) \times H^{\text{div}}(\Omega, \mathbb{R}^2)$. Consider the space of boundary port variables:
		\begin{equation}
		\begin{gathered}
		\mathcal{Z} := \left\{ \bm{z} \; \vert \; \bm{z} = \begin{pmatrix} \bm{f}_{\partial} \\ \bm{e}_{\partial} \end{pmatrix} \right\}, \quad \text{with} \quad
		\bm{f}_\partial = 
		\begin{pmatrix}
		w_t \\ \omega_{n} \\ \omega_{s} \\
		\end{pmatrix} \quad \text{and} \quad
		\bm{e}_\partial = 
		\begin{pmatrix}
		q_n \\ M_{nn} \\ M_{ns} \\
		\end{pmatrix}
		\end{gathered}
		\end{equation}
		The duality pairing between elements of $\mathcal{B}$ is then defined as follows:
		\begin{equation}
		\label{eq:bil_Min}
		\left\langle \left\langle (\bm{f}_1, \bm{e}_1, \bm{z}_1), (\bm{f}_2, \bm{e}_2, \bm{z}_2) \right\rangle \right\rangle :=  \left\langle \bm{e}_1, \bm{f}_2 \right\rangle_{\mathscr{L}^2(\Omega)}  +  \left\langle \bm{e}_2,  \bm{f}_1  \right\rangle_{\mathscr{L}^2(\Omega)} + \int_{\partial \Omega} B_J(\bm{z}_1, \bm{z}_2) \d{s} , 
		\end{equation}
		where the pairing $\left\langle \cdot, \cdot \right\rangle_{\mathscr{L}^2(\Omega)}$ is the $L^2$ inner product on space $\mathscr{L}^2(\Omega)$ and $B_J(\bm{z}_1, \bm{z}_2) := (\bm{f}_{\partial, 1})^T \bm{e}_{\partial, 2} + (\bm{f}_{\partial, 2}) ^T \bm{e}_{\partial, 1}$. 
		\begin{theorem}[Stokes-Dirac Structure for the Mindlin plate in tensorial form]
			Consider the space of power variables $\mathcal{B}$ defined in \eqref{eq:bondMin} and the matrix differential operator $J$ in~\eqref{eq:PH_sys_Min_Ten}. By theorem \ref{th:StokesDirac}, the linear subspace $\mathcal{D} \subset \mathcal{B}$:
			\begin{equation}
			\mathcal{D} =  \left\{(\bm{f}, \bm{e},\bm{z}) \in \mathcal{B} | \; \bm{f} = - \diffp{\bm{\alpha}}{t} = -J\bm{e}, \; \bm{z} = \begin{pmatrix} \bm{f}_{\partial} \\ \bm{e}_{\partial} \end{pmatrix} 
			\right\},
			\end{equation}
			is a Stokes-Dirac structure with respect to the pairing $\left\langle \left\langle \cdot, \cdot \right\rangle \right\rangle$ given by \eqref{eq:bil_Min}.
		\end{theorem}
		\begin{remark}
			The Hamiltonian formulation and associated Stokes-Dirac structure involve a Hamiltonian skew-symmetric differential operator $J$ on the space $\mathcal{E}$ of the coenergy variables. This represents a significant difference with respect to the Hamiltonian system obtained by the Legendre transformation of the Euler-Lagrange system \cite{jetMin}. In the former formulation the total dimension of this system is 8, in the latter 6. Furthermore, the latter system is ruled by an algebraic operator, the differential part being hidden in the variational derivative of the Hamiltonian, leading to the following system:
			\begin{equation}
			\begin{pmatrix}
			w\\ \theta_x\\ \theta_y\\p_w\\p_{\theta, x}\\p_{\theta, y}\\
			\end{pmatrix} = 
			\begin{bmatrix}
			0&0&0&1&0&0 \\
			0&0&0&0&1&0 \\
			0&0&0&0&0&1 \\
			-1&0&0&0&0&0\\
			0&-1&0&0&0&0\\
			0&0&-1&0&0&0\\
			\end{bmatrix}
			\begin{pmatrix}
			\delta_w H\\ \delta_{\theta, x} H\\ \delta_{\theta, y} H\\ \delta_{p_w} H\\\delta_{p_{\theta, x}} H\\\delta_{p_{\theta, y}} H\\
			\end{pmatrix},
			\end{equation}
		\end{remark}
		where $p_w$ is the linear momentum and $p_{\theta, x}, p_{\theta, y}$ are the angular momenta.
	}

	\section{Discretization of the Mindlin plate using a Partioned Finite Element Method}
	The Partioned Finite Element Method (PFEM) \cite{CardosoRibeiro2018} consists of putting the system into weak form first and of applying the integration by parts on a subset of the overall system second. For the Mindlin plate the integration by parts can be applied to first two lines of system \eqref{eq:PH_sys_Min_Ten}. This choice will make appear the momenta and forces at the boundary as control inputs. Alternatively, the last two of \eqref{eq:PH_sys_Min_Ten} could have been selected to perform the integration by parts. In this latter case the linear and angular velocities at the boundary would appear as control inputs. Keeping this in mind, the most suitable choice will depend on the physical problem under consideration.
	
	\subsection{Weak form}
	The test functions are of scalar, vectorial and tensorial nature. Keeping the same notation as in section \ref{sec:PH_ten_Min}, the scalar test function is denoted by $v_w$, the vectorial one by $\bm{v}_{\theta}$, $\bm{v}_{\gamma}$  the tensorial one by $\mathbb{V}_{\kappa}$.
	
	\subsubsection{Boundary control through forces and momenta}
	The first line of \eqref{eq:PH_sys_Min_Ten} is multiplied  by $v_w$ (multiplication by a scalar), the second line and the fourth by $\bm{v}_{\theta}$, $\bm{v}_{\gamma}$ (scalar product of $\mathbb{R}^2$) and the third one by $\mathbb{V}_{\kappa}$ (tensor contraction):
	\begin{align}
	\int_{\Omega} v_w \diffp{\alpha_w}{t}  \d\Omega &=  \int_{\Omega} v_w \mathrm{div}(\bm{e}_{\gamma}) \, \d\Omega,  \label{eq:wf1_min_ten}\\
	\int_{\Omega} \bm{v}_{\theta} \cdot \diffp{\bm\alpha_{\theta}}{t}   \d\Omega &= \int_{\Omega} \bm{v}_{\theta} \cdot (\mathrm{Div}( \mathbb{E}_{\kappa}) + \bm{e}_{\gamma}) \,  \d\Omega,  \label{eq:wf2_min_ten} \\
	\int_{\Omega} \mathbb{V}_{\kappa} \cddot \diffp{\mathbb{A}_{\kappa}}{t}   \d\Omega &= \int_{\Omega} \mathbb{V}_{\kappa} \cddot \mathrm{Grad}(\bm{e}_\theta) \, \d\Omega, \label{eq:wf3_min_ten} \\
	\int_{\Omega} \bm{v}_{\gamma} \cdot \diffp{\bm\alpha_{\gamma}}{t}   \d\Omega &= \int_{\Omega} \bm{v}_{\gamma} \cdot (\mathrm{grad}({e}_{w}) - \bm{e}_{\theta}) \, \d\Omega.  \label{eq:wf4_min_ten}
	\end{align}
	
	The right-hand side of equation \eqref{eq:wf1_min_ten} has to be integrated by parts
	\begin{equation}
	\label{eq:line1_ten_min}
	\int_{\Omega} v_w \mathrm{div}(\bm{e}_{\gamma})  \d\Omega = \int_{\partial \Omega} v_w \underbrace{\bm{n} \cdot \bm{e}_{\gamma}}_{q_n}  \d{s} - \int_{\Omega} \mathrm{grad}(v_w)  \cdot \bm{e}_{\gamma}  \d\Omega,\end{equation}
	as well as the right-hand side of equation \eqref{eq:wf2_min_ten}
	\begin{equation}
	\label{eq:line2_ten_min}
	\int_{\Omega} \bm{v}_{\theta} \cdot (\mathrm{Div}(\mathbb{E}_{\kappa}) + \bm{e}_{\gamma})  \d\Omega = \int_{\partial \Omega} \bm{v}_{\theta} \cdot (\bm{n} \cdot \mathbb{E}_{\kappa})  \d{s} -\int_{\Omega} \left\{ \mathrm{Grad}(\bm{v}_{\theta}) \cddot \mathbb{E}_{\kappa} - \bm{v}_{\theta} \cdot \bm{e}_{\gamma}\right\}  \d\Omega.
	\end{equation}
	
	The usual additional manipulation is performed on the boundary term containing the momenta, so that the proper boundary values arise:
	\begin{equation}
	\label{eq:line1_ipbc_ten_min}
	\begin{aligned}
	\int_{\partial \Omega} \bm{v}_{\theta} \cdot (\bm{n} \cdot \mathbb{E}_{\kappa})  \d{s} &= \int_{\partial \Omega} \left\{(\bm{v}_{\theta} \cdot \bm{n}) \bm{n} + (\bm{v}_{\theta} \cdot \bm{s}) \bm{s} \right\} \cdot (\bm{n} \cdot \mathbb{E}_{\kappa}) \,  \d{s} \\
	&= \int_{\partial \Omega} \left\{ v_{\omega_n} M_{nn} + v_{\omega_s} M_{ns} \right\}  \d{s}.
	\end{aligned}
	\end{equation}
	So defining $v_{\omega_n} :=\bm{v}_{\theta} \cdot \bm{n}$ and $v_{\omega_s} := \bm{v}_{\theta} \cdot \bm{s}$ the final weak form obtained from system \eqref{eq:PH_sys_Min_Ten} reads:
	\begin{equation}
	\label{eq:WF_Min_Dyn}
	\begin{cases}
	\displaystyle\int_{\Omega} v_w \diffp{\alpha_w}{t}  \d\Omega  &= - \displaystyle\int_{\Omega} \mathrm{grad}(v_w)  \cdot \bm{e}_{\gamma}  \d\Omega +  \displaystyle\int_{\partial \Omega} v_w q_n  \d{s}, \vspace{2mm}\\
	\displaystyle\int_{\Omega} \bm{v}_{\theta} \cdot \diffp{\bm\alpha_{\theta}}{t}   \d\Omega &=  -\displaystyle\int_{\Omega} \left\{ \mathrm{Grad}(\bm{v}_{\theta}) \cddot \mathbb{E}_{\kappa} - \bm{v}_{\theta} \cdot \bm{e}_{\gamma}\right\}  \d\Omega + \displaystyle\int_{\partial \Omega} \left\{ v_{\omega_n} M_{nn} + v_{\omega_s} M_{ns} \right\}  \d{s}, \vspace{2mm} \\
	\displaystyle\int_{\Omega} \mathbb{V}_{\kappa} \cddot \diffp{\mathbb{A}_{\kappa}}{t}   \d\Omega &= \displaystyle\int_{\Omega} \mathbb{V}_{\kappa} \cddot \mathrm{Grad}(\bm{e}_\theta)  \d\Omega,  \vspace{2mm} \\
	\displaystyle\int_{\Omega} \bm{v}_{\gamma} \cdot \diffp{\bm\alpha_{\gamma}}{t}   \d\Omega &= \displaystyle\int_{\Omega} \bm{v}_{\gamma} \cdot (\mathrm{grad}({e}_{w}) - \bm{e}_{\theta})  \d\Omega.
	\end{cases}
	\end{equation}
	In this first case,  the boundary controls $\bm{u}_\partial$ and the corresponding output $\bm{y}_\partial$ are: 
	\[\bm{u}_\partial = 
	\begin{pmatrix}
	q_n \\
	M_{nn} \\
	M_{ns} \\
	\end{pmatrix}_{\partial \Omega}, \qquad
	\bm{y}_\partial = 
	\begin{pmatrix}
	w_t \\
	\omega_n \\
	\omega_s \\
	\end{pmatrix}_{\partial \Omega}.
	\]

	\subsubsection{Boundary control through kinematic variables}
	Alternatively, in this second case, the same procedure can be performed on the two last lines of the system written in weak form (equations \eqref{eq:wf3_min_ten}, \eqref{eq:wf4_min_ten}). Once the due calculations are carried out, it is found:
	\begin{equation}
	\label{eq:WF_Min_Kin}
	\begin{cases}
	\displaystyle\int_{\Omega} v_w \diffp{\alpha_w}{t}  \d\Omega  &= \displaystyle \int_{\Omega} v_w \mathrm{div}(\bm{e}_{\gamma})  \d\Omega, \vspace{2mm}\\
	\displaystyle\int_{\Omega} \bm{v}_{\theta} \cdot \diffp{\bm\alpha_{\theta}}{t}   \d\Omega &= \displaystyle \int_{\Omega} \bm{v}_{\theta} \cdot (\mathrm{Div}(\mathbb{E}_{\kappa}) + \bm{e}_{\gamma}) \;  \d\Omega, \vspace{2mm} \\
	\displaystyle\int_{\Omega} \mathbb{V}_{\kappa} \cddot \diffp{\mathbb{A}_{\kappa}}{t}   \d\Omega &= - \displaystyle\int_{\Omega} \mathrm{Div}(\mathbb{V}_{\kappa}) \cdot \bm{e}_\theta \;  \d\Omega +  \displaystyle\int_{\partial \Omega} \left\{ v_{M_{nn}} \omega_{n} + v_{M_{ns}} \omega_{s} \right\}  \d{s}, \vspace{2mm} \\
	\displaystyle\int_{\Omega} \bm{v}_{\gamma} \cdot \diffp{\bm\alpha_{\gamma}}{t}   \d\Omega &=  - \displaystyle\int_{\Omega} \left\{ \mathrm{div}(\bm{v}_{\gamma}) e_w + \bm{v}_{\gamma} \cdot \bm{e}_{\theta} \right\} \, \d\Omega + \displaystyle\int_{\partial \Omega} v_{q_{n}} w_t \;  \d{s},
	\end{cases}
	\end{equation}
	
	where $v_{M_{nn}} = \mathbb{V}_{\kappa} \cddot (\bm{n} \otimes \bm{n}), \;  v_{M_{ns}} = \mathbb{V}_{\kappa} \cddot (\bm{s} \otimes \bm{n}) \;$ and $v_{q_n} = \bm{v}_{\gamma} \cdot \bm{n}$.
	In this second case,  the boundary controls $\bm{u}_\partial$ and corresponding output $\bm{y}_\partial$ are:
	\[\bm{u}_\partial = 
	\begin{pmatrix}
	{w_t} \\
	\omega_{n} \\
	\omega_{s} \\
	\end{pmatrix}_{\partial \Omega}, \qquad
	\bm{y}_\partial = 
	\begin{pmatrix}
	q_n \\
	M_{nn} \\
	M_{ns} \\
	\end{pmatrix}_{\partial \Omega}.
	\]
	
	\subsection{Finite-dimensional port-Hamiltonian system}
	\label{sec:FD_system}
	{
		In this section the formulation \eqref{eq:WF_Min_Dyn} is used is order to explain the discretization procedure. Test and co-energy variables are discretized using the same basis functions (Galerkin Method):
		\begin{equation}
		\begin{aligned}
		v_w &= \sum_{i = 1}^{N_w} \phi_w^i(x,y) \, v_w^i, \\
		\bm{v}_\theta &= \sum_{i = 1}^{N_\theta} \bm\phi_\theta^i(x,y) \, v_\theta^i, \\
		\mathbb{V}_\kappa &= \sum_{i = 1}^{N_\kappa} \bm\Phi_\kappa^i(x,y) \, v_\kappa^i,\\
		\bm{v}_{\gamma} &= \sum_{i = 1}^{N_{\gamma}} \bm\phi_{\gamma}^i(x,y) \, v_{\gamma}^i,\\
		\end{aligned} \qquad \quad
		\begin{aligned}
		e_w &= \sum_{i = 1}^{N_w} \phi_w^i(x,y) \, e_w^i(t),  \\
		\bm{e}_\theta &= \sum_{i = 1}^{N_\theta} \bm\phi_\theta^i(x,y) \, e_\theta^i(t), \\
		\mathbb{E}_\kappa &= \sum_{i = 1}^{N_\kappa} \bm\Phi_\kappa^i(x,y) \, e_\kappa^i(t),\\
		\bm{e}_{\gamma} &= \sum_{i = 1}^{N_{\gamma}} \bm\phi_{\gamma}^i(x,y) \, e_{\gamma}^i(t).\\
		\end{aligned} 
		\end{equation}
		The basis functions $\phi_w^i, \, \bm\phi_\theta^i, \, \bm\Phi_\kappa^i, \, \bm\phi_{\gamma}^i$ have to be chosen in a suitable function space $\mathcal{V}^h$ in the domain of operator $J$, defined in \eqref{eq:PH_sys_Min_Ten}, i.e. $\mathcal{V}^h \subset \mathcal{V} \in \mathcal{D}(J)$. This will be discussed in section \ref{sec:Num}. The discretized skew-symmetric bilinear form on the right-hand side of \eqref{eq:WF_Min_Dyn} then becomes:
		\begin{equation}
		\bm{J}_d = 
		\begin{bmatrix}
		0 & 0 & 0 & -\bm{D}_{\mathrm{grad}}^T      \vspace{.3mm}\\ 
		0 & 0 & -\bm{D}_{\mathrm{Grad}}^T & -\bm{D}_0^T \vspace{.3mm}\\
		0 & \bm{D}_{\mathrm{Grad}} & 0 & 0         \vspace{.3mm}\\
		\bm{D}_{\mathrm{grad}} & \bm{D}_0 & 0 & 0       \vspace{.3mm}\\
		\end{bmatrix},
		\end{equation}
		where the matrices are computed in the following way:
		\begin{equation}
		\begin{aligned}
		\bm{D}_{\mathrm{Grad}}(i,j) &= \int_{\Omega} \bm{\Phi}_{\kappa}^i : \mathrm{Grad} (\bm{\phi}_{\theta}^j) \d\Omega, \quad \in \mathbb{R}^{N_\kappa \times N_\theta},\\
		\bm{D}_{\mathrm{grad}}(i,j) &= \int_{\Omega} \bm{\phi}_{\gamma}^i \cdot  \mathrm{grad}(\bm{\phi}^j_{w}) \d\Omega, \quad \in \mathbb{R}^{N_{\gamma} \times N_w},\\
		\bm{D}_0(i,j)  &= -\int_{\Omega} \bm{\phi}_{\gamma}^i \cdot  \bm{\phi}^j_{\theta} \d\Omega, \quad \in \mathbb{R}^{N_{\gamma} \times N_\theta}. \\
		\end{aligned}
		\end{equation}
		The notation $\bm{D}(i,j)$ indicates the entry in matrix $\bm{D}$ corresponding to the $i \, {\text{th}}$ row and $j \,{\text{th}}$ column. The energy variables are deduced from the co-energy variables:
		\begin{equation}
		\begin{aligned}
		\alpha_w &= \rho h e_w, \\
		\bm\alpha_{\theta} &=  \frac{\rho h^3}{12} \bm{e}_{\theta}, \\
		\end{aligned} \qquad
		\begin{aligned}
		\mathbb{A}_{\kappa} &= \mathbb{D}^{-1} \mathbb{E}_{\kappa}, \\
		\bm\alpha_{\gamma} &= \frac{1}{G h k} \bm{e}_{\gamma}. \\
		\end{aligned}
		\end{equation}
		The symmetric bilinear form on the left-hand side of \eqref{eq:WF_Min_Dyn} is discretized as:
		\begin{equation}
		\begin{gathered}
		\bm{M} = \text{diag}[\bm{M}_w,\, \bm{M}_\theta,\, \bm{M}_\kappa,\, \bm{M}_{\gamma}], \quad \text{with} \\
		\begin{aligned}
		&\bm{M}_w(i,j) = \int_{\Omega} \rho h \, \bm{\phi}_w^i \, \bm{\phi}_w^j \d\Omega, \; \in \mathbb{R}^{N_w \times N_w}, \\
		&\bm{M}_\theta(i,j) = \int_{\Omega} \frac{\rho h^3}{12} \bm{\phi}_\theta^i \cdot \bm{\phi}_\theta^j \d\Omega, \; \in \mathbb{R}^{N_\theta \times N_\theta},\\
		\end{aligned} \qquad
		\begin{aligned}
		&\bm{M}_\kappa(i,j) = \int_{\Omega}  \left( \mathbb{D}^{-1} \bm{\Phi}_\kappa^i \right) \cddot \bm{\Phi}_\kappa^j \d\Omega, \; \in \mathbb{R}^{N_\kappa \times N_\kappa},\\
		&\bm{M}_{\gamma}(i,j) = \int_{\Omega} \frac{1}{G h k} \bm{\phi}_{\gamma}^i \cdot \bm{\phi}_{\gamma}^j \d\Omega, \; \in \mathbb{R}^{N_{\gamma} \times N_{\gamma}}. \\
		\end{aligned}
		\end{gathered}
		\end{equation}
		The boundary variables are then discretized as:
		\begin{equation}
		\begin{aligned}
		q_n = \sum_{i = 1}^{N_{q_n}} \phi_{q_n}^i(s) q_n^i, \qquad
		M_{nn} = \sum_{i = 1}^{N_{M_{nn}}} \phi_{M_{nn}}^i(s) M_{nn}^i, \qquad
		M_{ns} = \sum_{i = 1}^{N_{M_{ns}}} \phi_{M_{ns}}^i(s) M_{ns}^i.
		\end{aligned}
		\end{equation}
		The variables are defined only over the boundary $\partial\Omega$. Consequently, the input matrix reads:
		\begin{equation}
		\bm{B} = \begin{bmatrix}
		\bm{B}_{q_n} & 0 & 0 \\
		0 & \bm{B}_{M_{nn}} & \bm{B}_{M_{ns}} \\
		0 & 0 & 0 \\
		0 & 0 & 0 \\
		\end{bmatrix}.
		\end{equation}
		The inner components are computed as:
		\begin{equation}
		\begin{aligned}
		\bm{B}_{q_n}(i,j) &= \int_{\partial\Omega} {\phi}_{w}^i \, {\phi}_{q_n}^j \d{s}, \quad \in \mathbb{R}^{N_w \times N_{q_n}} , \\
		\bm{B}_{M_{nn}}(i,j) &= \int_{\partial\Omega} (\bm{\phi}_{\theta}^i \cdot \bm{n})\, \phi_{M_{nn}}^j \d{s}, \quad \in \mathbb{R}^{N_w \times N_{M_{nn}}} , \\
		\bm{B}_{M_{ns}}(i,j) &= \int_{\partial\Omega} (\bm{\phi}_{\theta}^i \cdot \bm{s})\, \phi_{M_{ns}}^j \d{s}, \quad \in \mathbb{R}^{N_w \times N_{M_{ns}}}. \\
		\end{aligned}
		\end{equation}
		The final port-Hamiltonian system (as defined in \cite{beattie2018linear}, where the presence of a mass matrix $\bm{M}$ is taken into account) is written as:
		\begin{equation}
		\label{eq:PHdiscr_Min}
		\begin{aligned}
		\bm{M} \dot{\bm{e}} &= \bm{J}_d  \,\bm{e} + \bm{B} \, \bm{u}_{\partial}, \\
		\bm{y}_{\partial} &= \bm{B}^T \, \bm{e},
		\end{aligned} 
		\end{equation}
		where $\bm{e} = \left(e_w^1, \cdots, e_{\gamma}^{N_{\gamma}}\right)^T$ and $\bm{u}_{\partial} = \left(q_n^1, \dots, M_{ns}^{N_{M_{ns}}}\right)^T$ are the concatenations of the degrees of freedom for the different variables. The discrete Hamiltonian is then found as:
		\begin{equation}
		\label{eq:Hd}
		\begin{aligned}
		H_d &= \frac{1}{2} \int_{\Omega} \left\{ \alpha_w e_w + \bm{\alpha}_\theta \cdot \bm{e}_{\theta} + \mathbb{A}_{\kappa} \cddot \mathbb{E}_\kappa + \bm\alpha_{\gamma} \cdot \bm{e}_{\gamma}  \right\} \d\Omega \\
		&=  \frac{1}{2}  \left\{ \bm{e}_w^T \, \bm{M}_w \, \bm{e}_w + \bm{e}_\theta^T \, \bm{M}_\theta \, \bm{e}_\theta + \bm{e}_\kappa^T \, \bm{M}_\kappa \, \bm{e}_\kappa + \bm{e}_{\gamma}^T \, \bm{M}_{\gamma} \, \bm{e}_{\gamma}  \right\} \\
		&=  \frac{1}{2} \, \bm{e}^T  \bm{M}  \bm{e}.
		\end{aligned}
		\end{equation}
	}
	Then it naturally follows that:
	\begin{equation}
	\dot{H}_d = \bm{y}_{\partial}^T \bm{u}_{\partial}.
	\end{equation}
	This is equivalent to the energy balance of the continuous system, expressed by equation \eqref{eq:PowRate}. Definition \eqref{eq:Hd}, together with system \eqref{eq:PHdiscr_Min} are the finite-dimensional equivalent of \eqref{eq:H_min} and  \eqref{eq:PH_sys_Min_Ten}. Again, the discretized system obtained via PFEM shares the same properties as those of the original infinite-dimensional system, the discretization method is therefore structure-preserving. \\

	\subsection{Mixed boundary conditions}
	\label{sec:mix_bd}
\begin{figure}[t]
		\centering
		\includegraphics[width=0.7\textwidth]{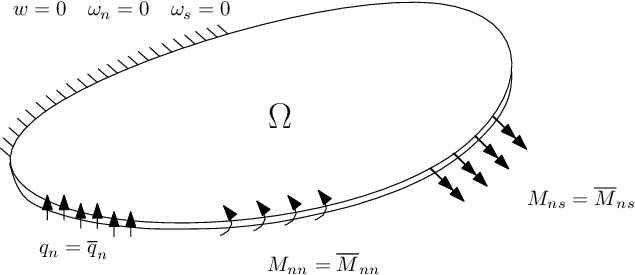}
		\caption{Clamped plate with distributed shear forces and torsional and flexural momenta at the border.}
		\label{fig:mix_bcs}
	\end{figure}

	This formulation provides as control term the dynamic variables at the boundary, namely forces and momenta, whereas the kinematic variables do not appear. In order to handle mixed boundary conditions, i.e. a clamped or a simply supported plate with free or loaded edges (see Fig. \ref{fig:mix_bcs}), the boundary control term has to be split into known and unknown variables, i.e. given boundary conditions and reactions at the boundaries respectively. These terms may be rearranged by introducing a permutation matrix $P$. The control term may be rewritten in the following way:
	\begin{equation}
	\bm{u}_{\partial} = \bm{B} \bm{P} 
	\begin{pmatrix}
	\bm{f} \\
	\bm{\lambda} \\
	\end{pmatrix} = 
	\bm{B} \begin{bmatrix}
	\bm{P}_{f} & \bm{P}_{\lambda} \\
	\end{bmatrix}
	\begin{pmatrix}
	\bm{f} \\
	\bm{\lambda} \\
	\end{pmatrix}.
	\end{equation}
	Equivalently the boundary outputs can be split. The terms corresponding to $\bm\lambda$ will be the kinematic variables set by the boundary conditions. The term corresponding to $\bm{f}$ are the kinematic variables at the controlled boundaries. The following equation allows splitting up the outputs into these two contributions:
	\begin{equation}
	\bm{y}_{\partial} = \bm{P}
	\begin{pmatrix}
	\bm{y}_{f} \\
	\bm{y}_{\lambda} \\
	\end{pmatrix} = 
	\begin{bmatrix}
	\bm{P}_{f} & \bm{P}_{\lambda} \\
	\end{bmatrix}
	\begin{pmatrix}
	\bm{y}_{f} \\
	\bm{y}_{\lambda} \\
	\end{pmatrix}.
	\end{equation}
	So that the output equation becomes:
	\begin{equation}
	\begin{pmatrix}
	\bm{y}_{f} \\
	\bm{y}_{\lambda} \\
	\end{pmatrix} = 
	\begin{bmatrix}
	\bm{P}_{f}^T \\
	\bm{P}_{\lambda}^T \\
	\end{bmatrix}
	\bm{B}^T \bm{e}.
	\end{equation}
	
	The port-Hamiltonian finite-dimensional system is rewritten equivalently by highlighting the known control terms, the Lagrange multipliers (reactions at the boundary) and the constraints, arising from the fact that {the function $\bm{y}_{\lambda}$ is known. In case of a homogenous Dirichlet condition, i.e. $\bm{y}_{\lambda} = 0$, it is obtained:
		\begin{equation}
		\label{eq:PHdiscr_mixed}
		\begin{aligned}
		\begin{bmatrix}
		\bm{M} & 0 \\
		0 & 0 \\
		\end{bmatrix}
		\frac{\d}{\d t}
		\begin{pmatrix}
		\bm{e}\\
		\bm{\lambda}\\
		\end{pmatrix} &=
		\begin{bmatrix}
		\bm{J}_d & \bm{B} \bm{P}_{\lambda} \\
		-\bm{P}_{\lambda}^T \bm{B}^T & 0 \\
		\end{bmatrix}
		\begin{pmatrix}
		\bm{e}\\
		\bm{\lambda}\\
		\end{pmatrix} + 
		\begin{bmatrix}
		\bm{B} \bm{P}_{f} \\
		0 \\
		\end{bmatrix} \bm{f}, \\
		\bm{y}_{f} &= \begin{bmatrix} \bm{P}_{f}^T \bm{B}^T \; 0 \end{bmatrix} \begin{pmatrix}
		\bm{e}\\
		\bm{\lambda}\\
		\end{pmatrix}. 
		\end{aligned}
		\end{equation}
	}
	This differential-algebraic port-Hamiltonian system can be treated applying results detailed in \cite{vanderSchaft2013,beattie2018linear}.
	
\begin{figure}[p]%
		\minipage{0.25\textwidth}%
		\includegraphics[width=\linewidth]{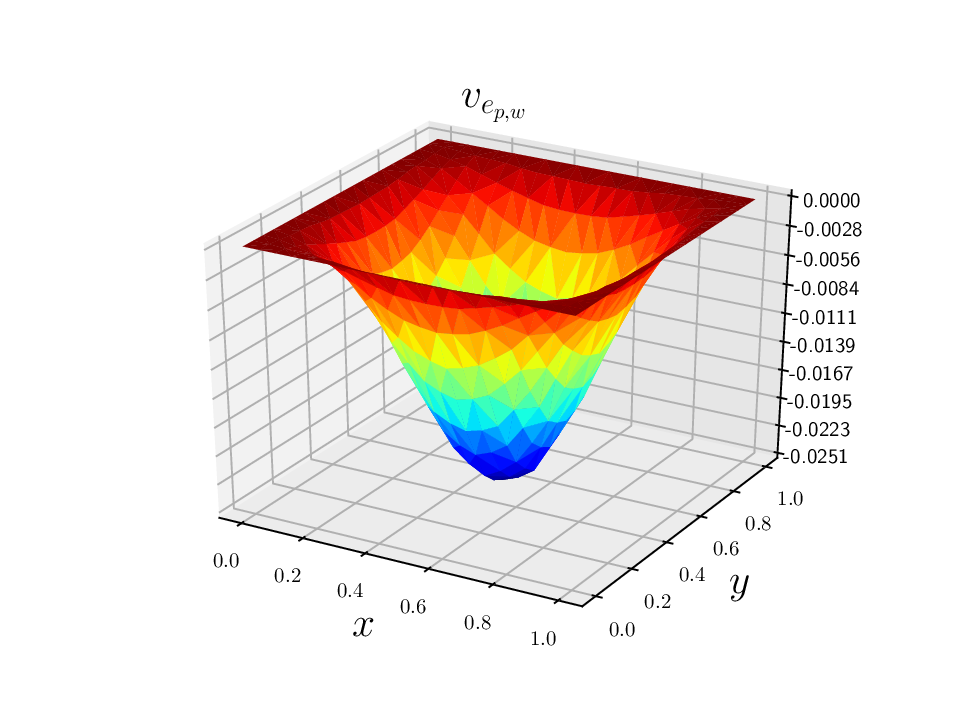}
		\caption*{$\widehat{\omega}_{11}$}\label{fig:CCCC1}%
		\endminipage
		\minipage{0.25\textwidth}%
		\includegraphics[width=\linewidth]{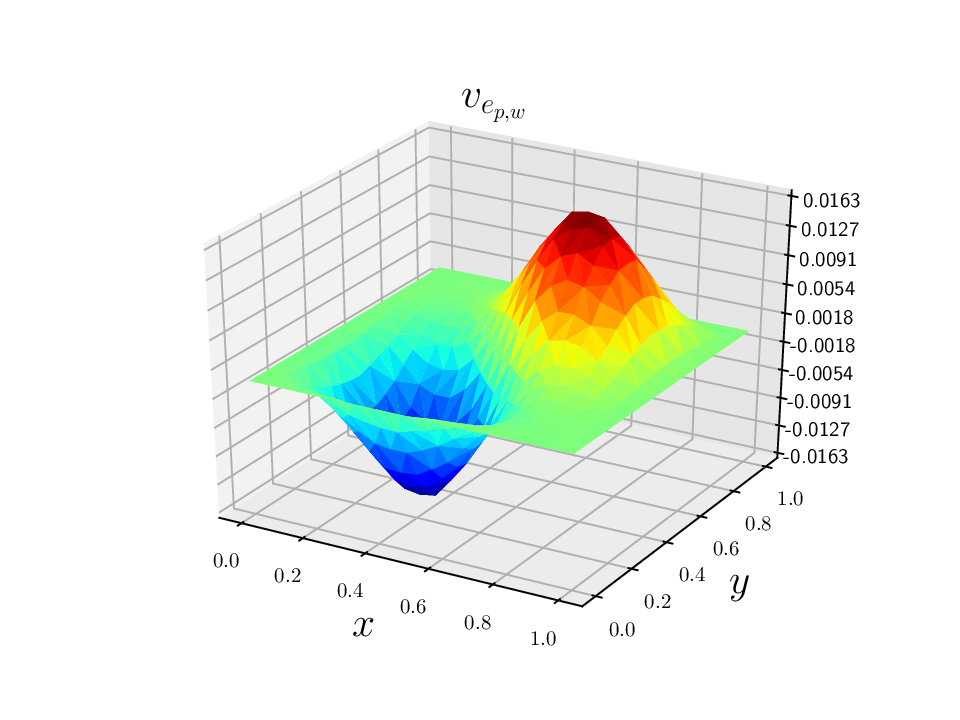}
		\caption*{$\widehat{\omega}_{21}$}\label{fig:CCCC2}%
		\endminipage
		\minipage{0.25\textwidth}%
		\includegraphics[width=\linewidth]{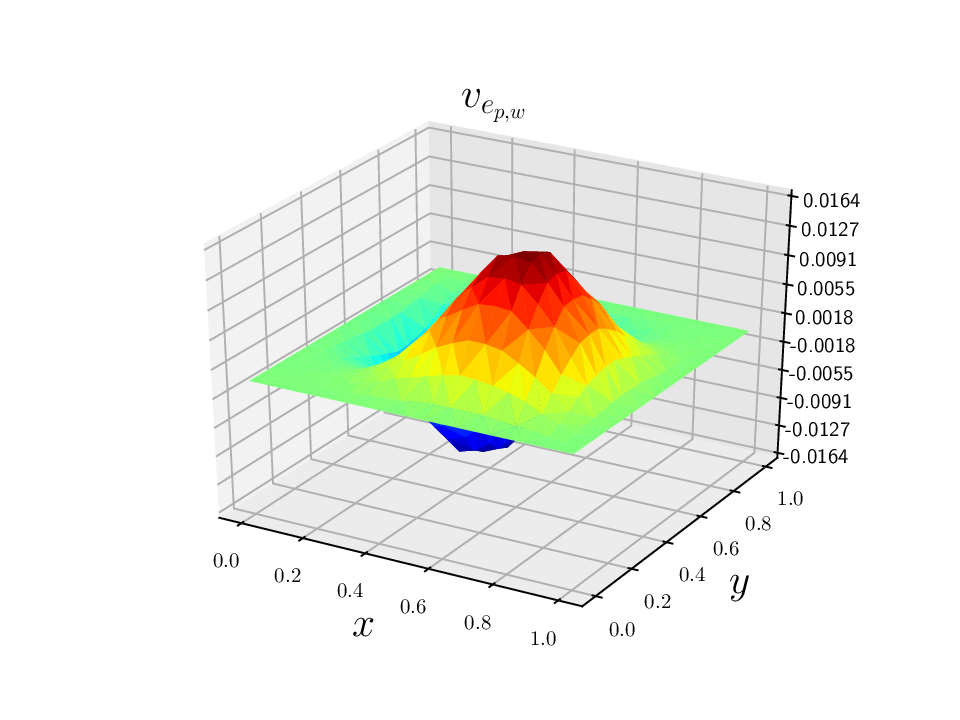}
		\caption*{$\widehat{\omega}_{12}$}\label{fig:CCCC3}%
		\endminipage
		\minipage{0.25\textwidth}%
		\includegraphics[width=\linewidth]{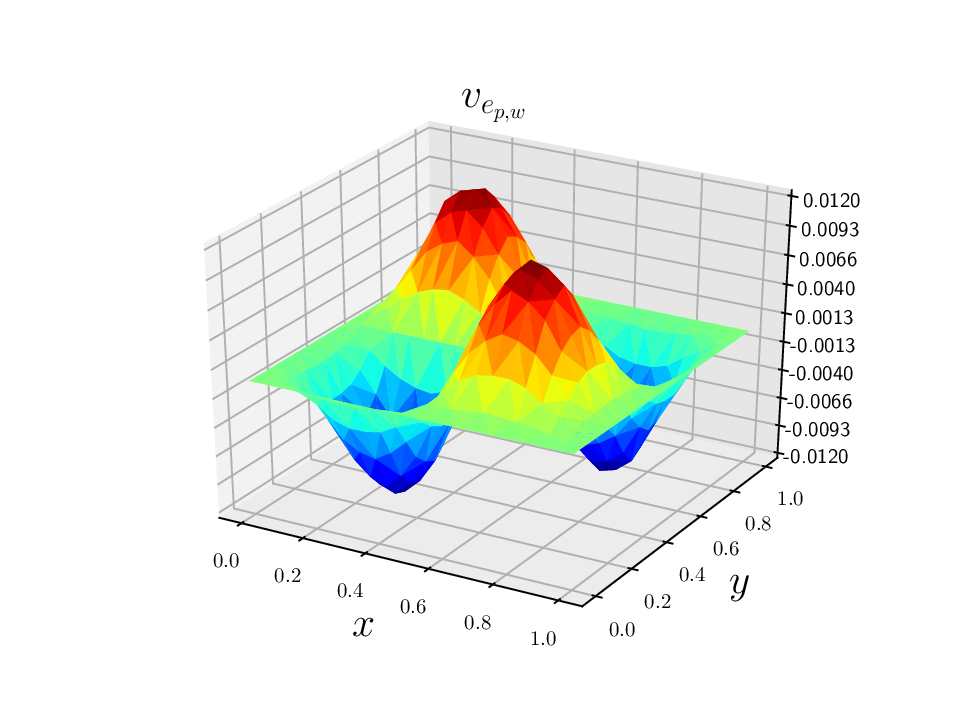}
		\caption*{$\widehat{\omega}_{22}$}\label{fig:CCCC4}%
		\endminipage
		\caption[Eigenvectors for CCCC]{Eigenvectors for the CCCC case.}%
		\label{fig:CCCC}%
	\end{figure}

\begin{figure}[p]%
		\minipage{0.25\textwidth}%
		\includegraphics[width=\linewidth]{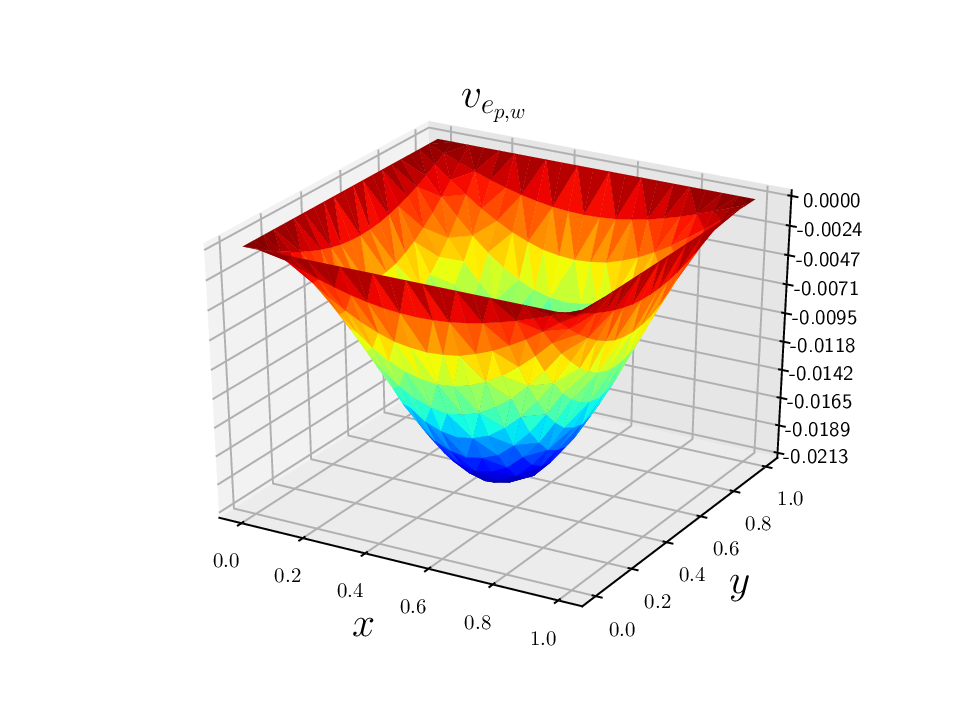}
		\caption*{$\widehat{\omega}_{11}$}\label{fig:SSSS1}%
		\endminipage
		\minipage{0.25\textwidth}%
		\includegraphics[width=\linewidth]{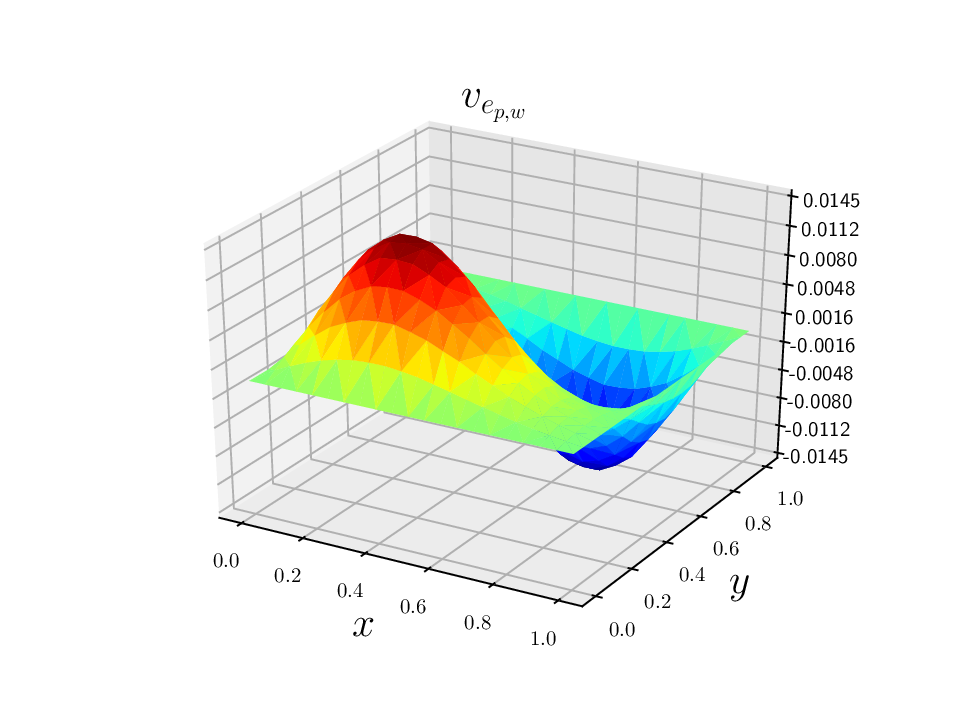}
		\caption*{$\widehat{\omega}_{21}$}\label{fig:SSSS2}%
		\endminipage
		\minipage{0.25\textwidth}%
		\includegraphics[width=\linewidth]{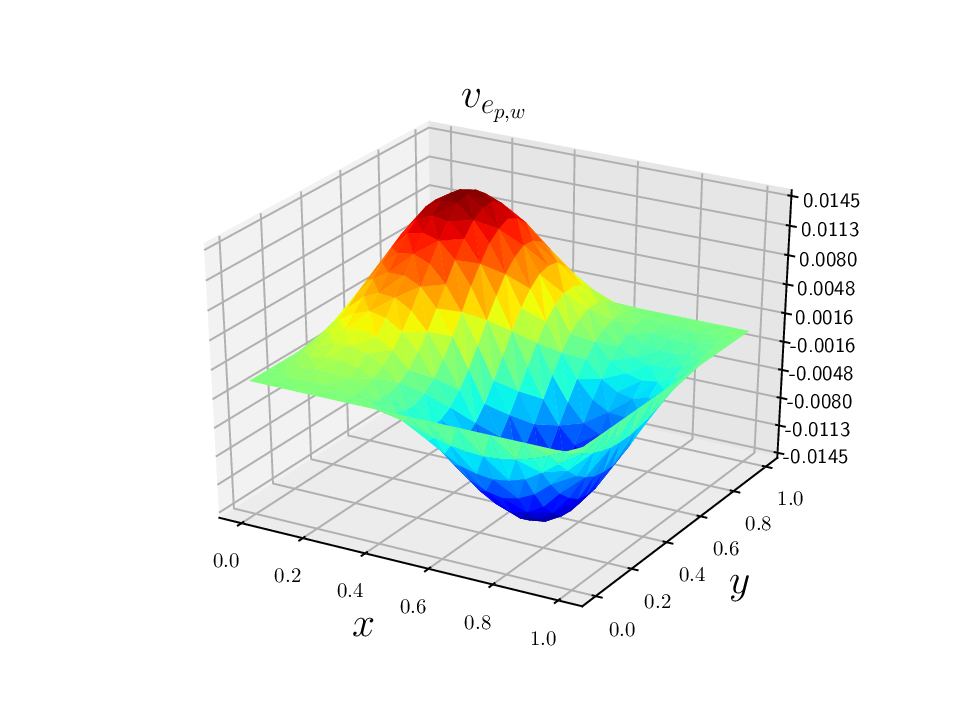}
		\caption*{$\widehat{\omega}_{12}$}\label{fig:SSSS3}%
		\endminipage 
		\minipage{0.25\textwidth}%
		\includegraphics[width=\linewidth]{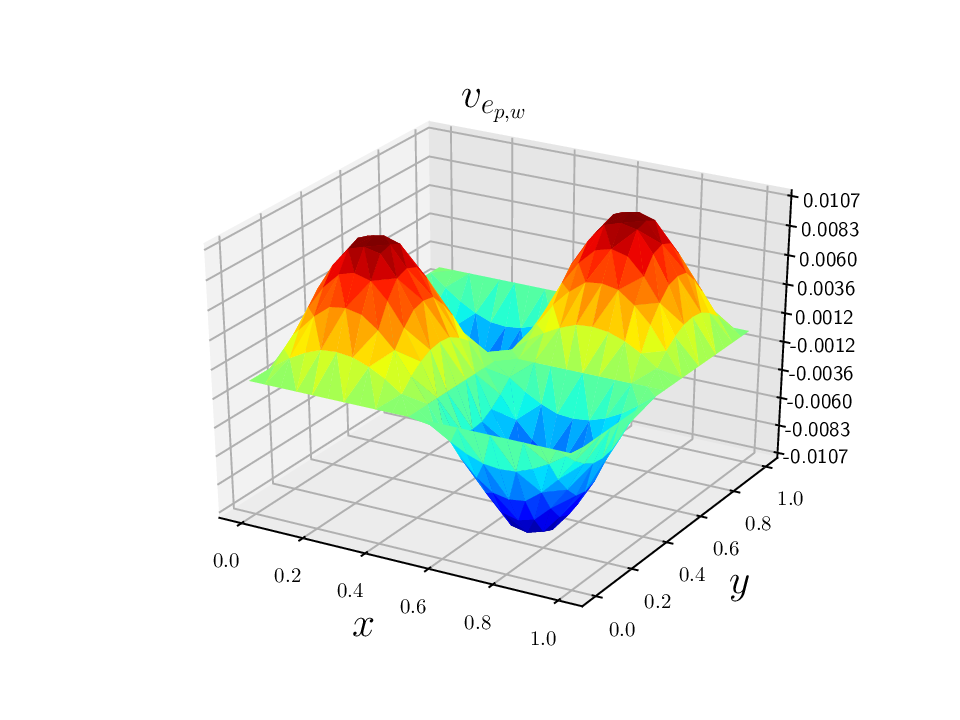}
		\caption*{$\widehat{\omega}_{22}$}\label{fig:SSSS4}%
		\endminipage
		\caption[Eigenvectors for CCCC]{Eigenvectors for the SSSS case.}%
		\label{fig:SSSS}%
	\end{figure}

\begin{figure}[p]%
		\minipage{0.25\textwidth}%
		\includegraphics[width=\linewidth]{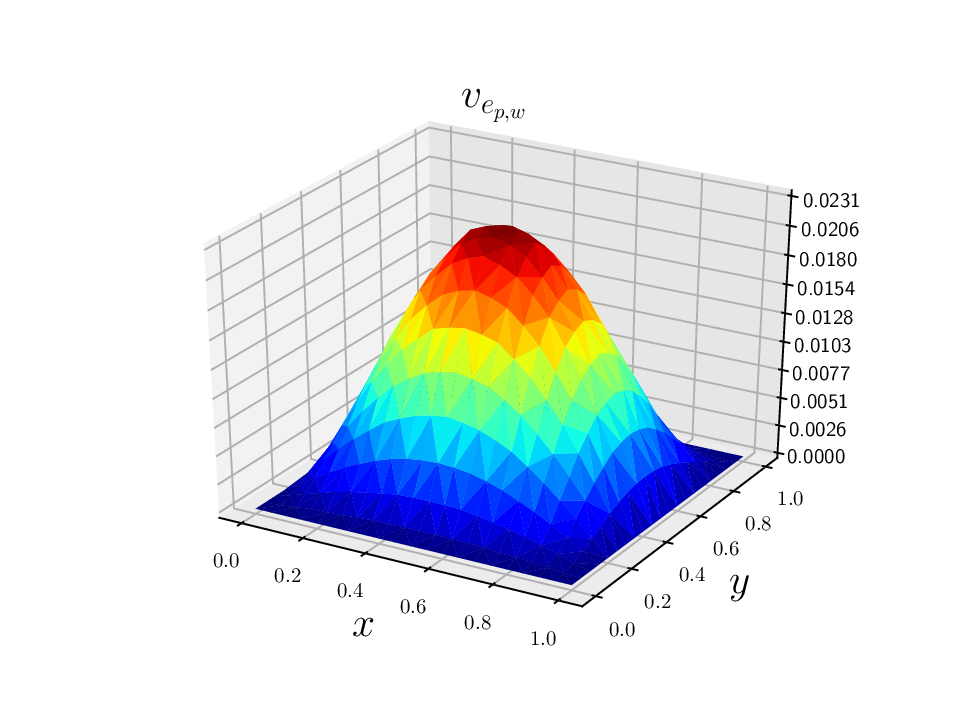}
		\caption*{$\widehat{\omega}_{11}$}\label{fig:SCSC1}%
		\endminipage
		\minipage{0.25\textwidth}%
		\includegraphics[width=\linewidth]{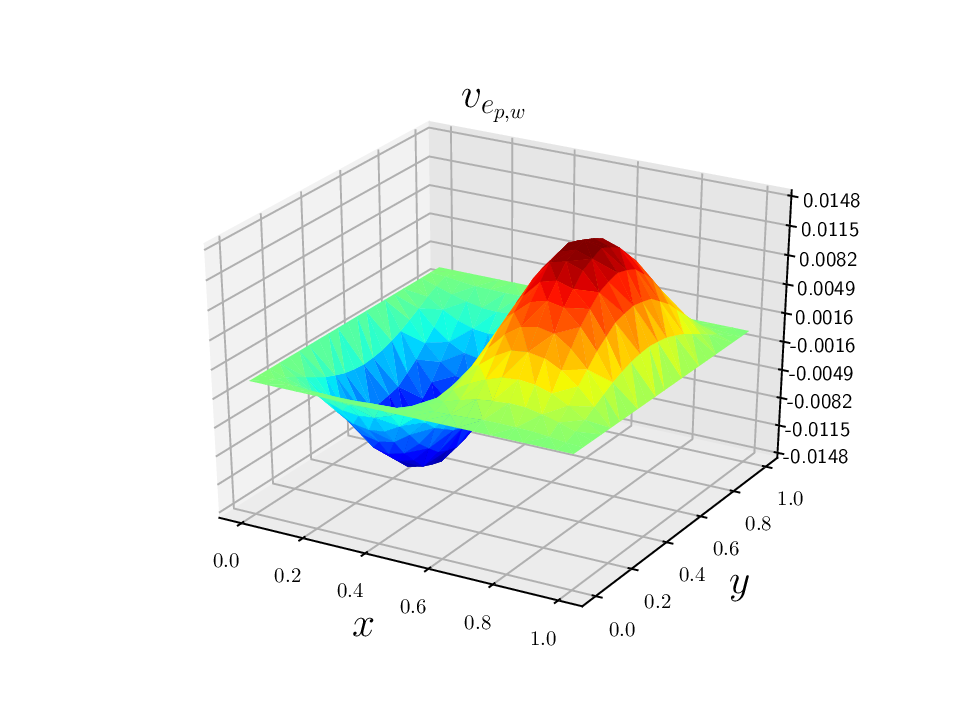}
		\caption*{$\widehat{\omega}_{21}$}\label{fig:SCSC2}%
		\endminipage
		\minipage{0.25\textwidth}%
		\includegraphics[width=\linewidth]{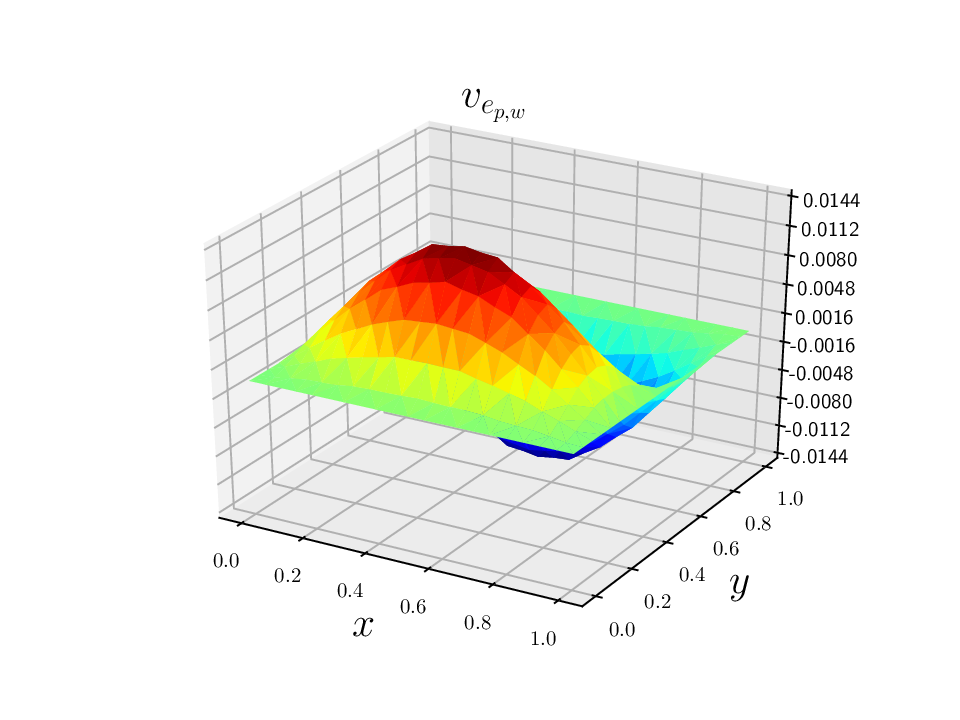}
		\caption*{$\widehat{\omega}_{12}$}\label{fig:SCSC3}%
		\endminipage 
		\minipage{0.25\textwidth}%
		\includegraphics[width=\linewidth]{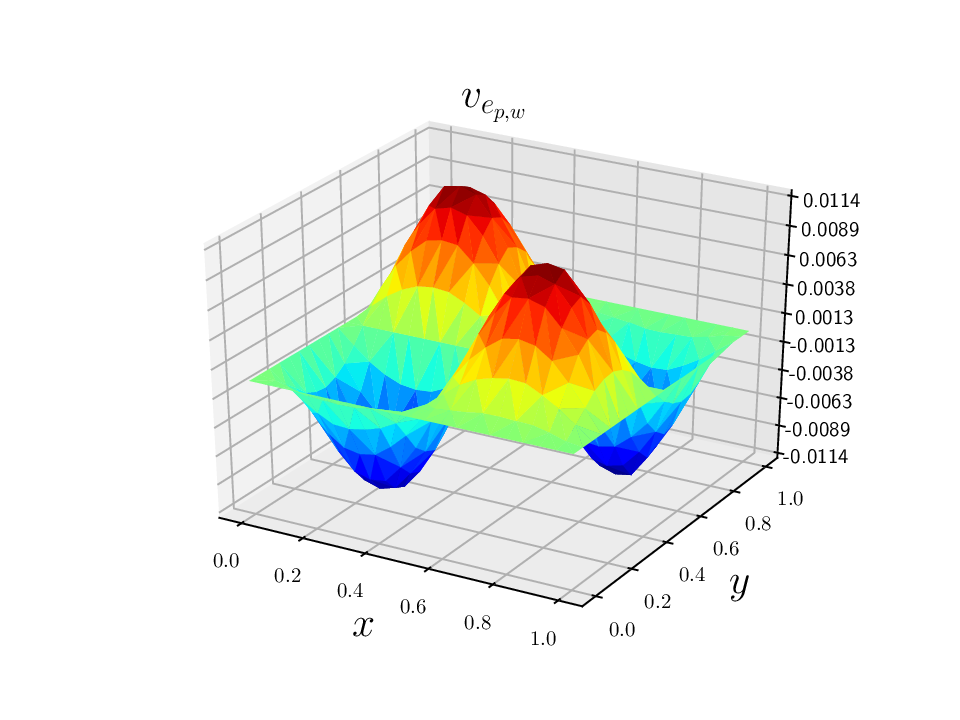}
		\caption*{$\widehat{\omega}_{22}$}\label{fig:SCSC4}%
		\endminipage
		\caption[Eigenvectors for CCCC]{Eigenvectors for the SCSC case.}%
		\label{fig:SCSC}%
	\end{figure}

\begin{figure}[p]%
		\minipage{0.25\textwidth}%
		\includegraphics[width=\linewidth]{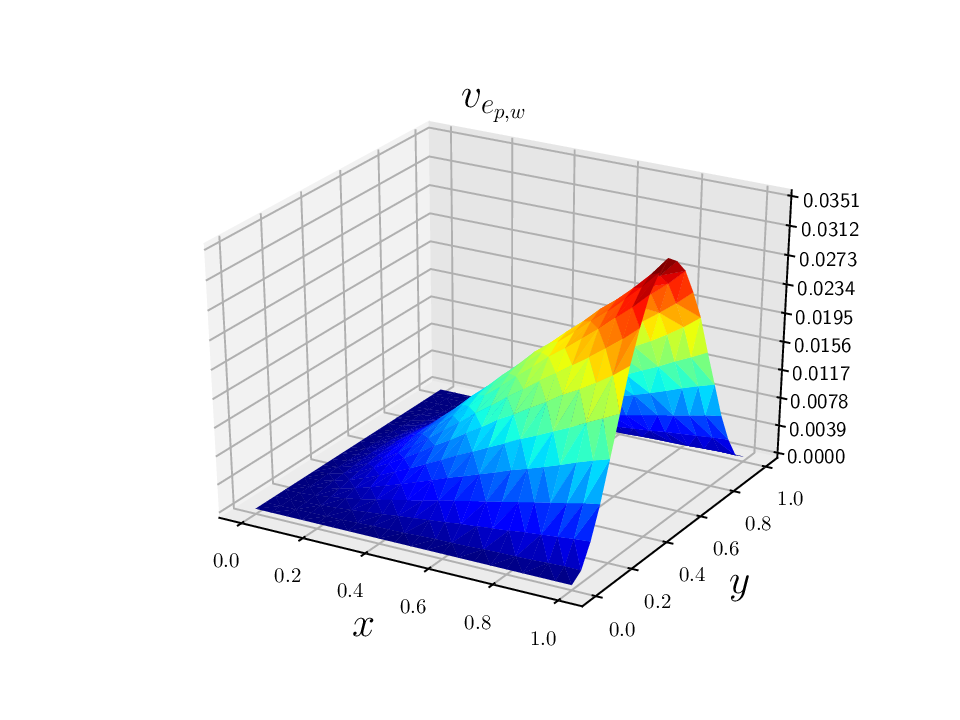}
		\caption*{$\widehat{\omega}_{\frac{1}{2}1}$}\label{fig:CCCF1}%
		\endminipage
		\minipage{0.25\textwidth}%
		\includegraphics[width=\linewidth]{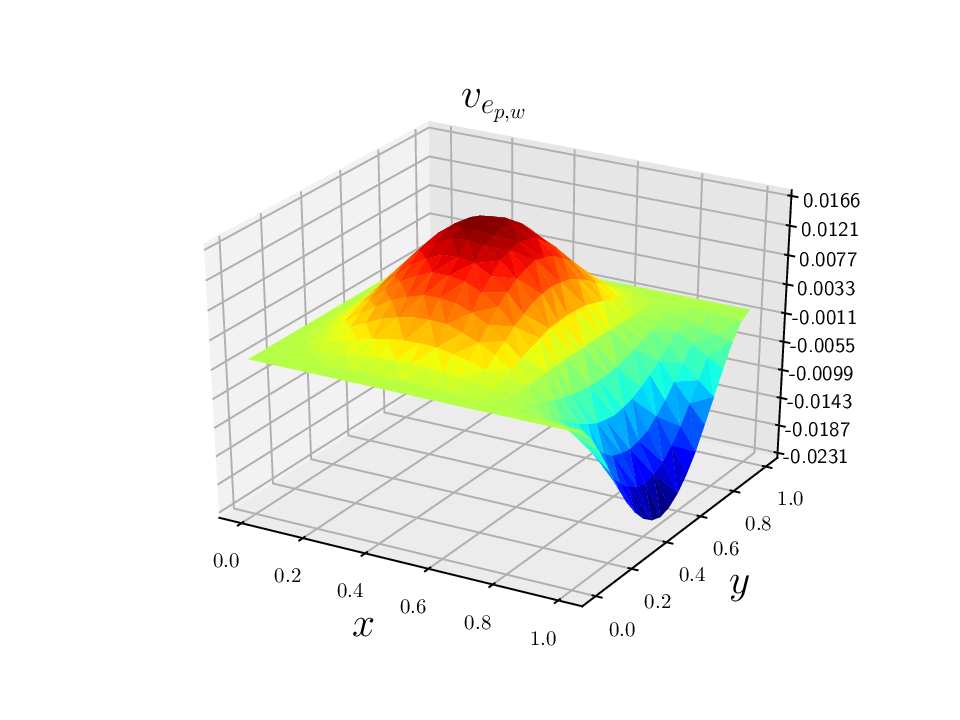}
		\caption*{$\widehat{\omega}_{\frac{3}{2}1}$}\label{fig:CCCF2}%
		\endminipage
		\minipage{0.25\textwidth}%
		\includegraphics[width=\linewidth]{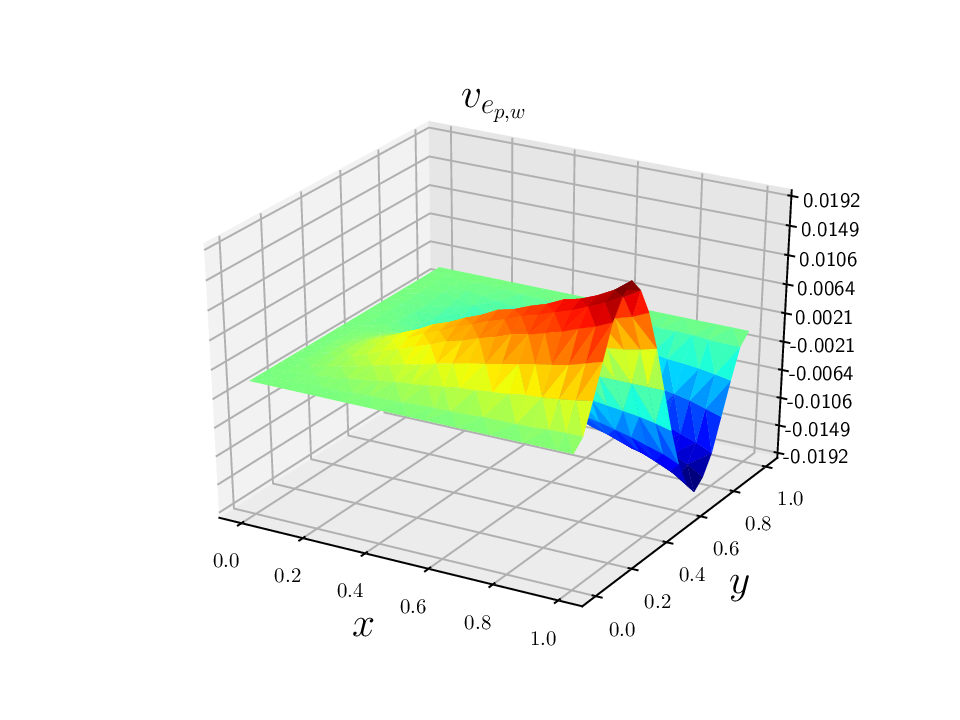}
		\caption*{$\widehat{\omega}_{\frac{1}{2}2}$}\label{fig:CCCF3}%
		\endminipage 
		\minipage{0.25\textwidth}%
		\includegraphics[width=\linewidth]{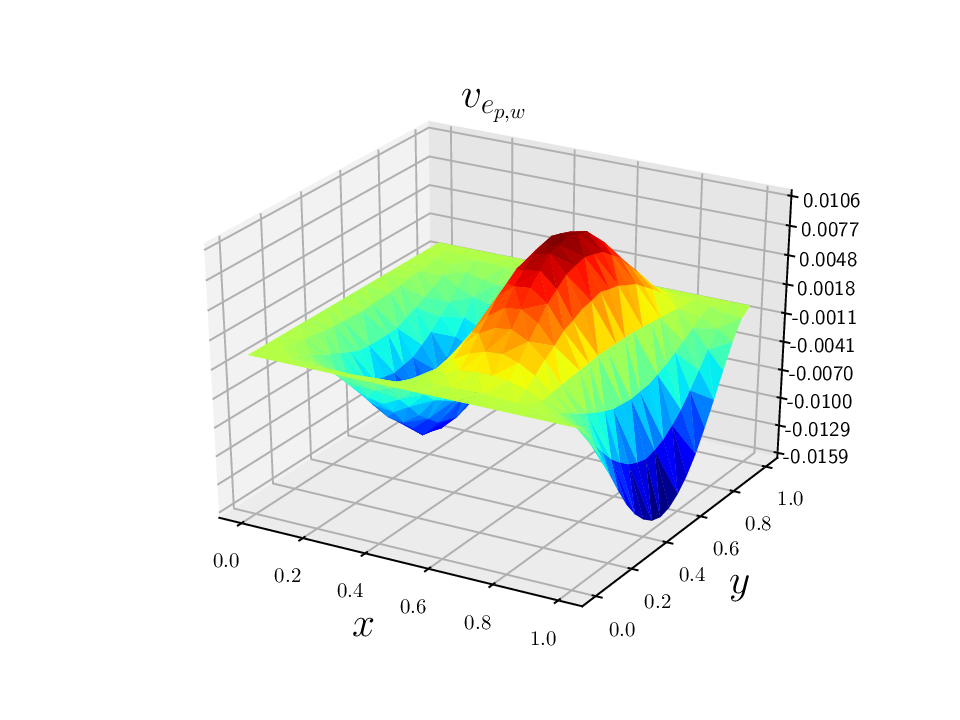}
		\caption*{$\widehat{\omega}_{\frac{5}{2}1}$}\label{fig:CCCF4}%
		\endminipage
		\caption[Eigenvectors for CCCC]{Eigenvectors for the CCCF case.}%
		\label{fig:CCCF}%
	\end{figure}

	{
		\section{Numerical studies}
		\label{sec:Num}
		In this section, the consistency of the proposed model is illustrated numerically. For this purpose the computation of the eigenvalues of a square plate and time-domain simulations for several boundary conditions are presented. The formulation used in the numerical implementation is the one presented in section \ref{sec:FD_system}. 
		\subsection{Finite Element Choice}
		The domain of the operator $J$ in \eqref{eq:PH_sys_Min_Ten} is 
		\begin{equation*}
		\label{eq:domJ}
		\mathcal{D}(J) = H^{1}(\Omega) \times H^{1}(\Omega, \mathbb{R}^2) \times H^{\text{Div}}(\Omega, \mathbb{R}^{2 \times 2}_{\text{sym}}) \times H^{\text{div}}(\Omega, \mathbb{R}^2) \quad \text{and boundary conditions}.
		\end{equation*}
		For this reason a suitable choice for the functional space is:
		\begin{equation}
		\label{eq:funcSpace}
		(v_w, \, \bm{v}_\theta, \, \mathbb{V}_\kappa, \, \bm{v}_{\gamma}) \in H^{1}(\Omega) \times H^{1}(\Omega, \mathbb{R}^2) \times H^{1}(\Omega, \mathbb{R}^{2 \times 2}_{\text{sym}}) \times H^{1}(\Omega, \mathbb{R}^2) \equiv \mathscr{H},
		\end{equation}
		since $\mathscr{H} \subset \mathcal{D}(J)$. Then, for the Finite Element choice, denote:
		\[ H_r^1(\mathbb{P}_l) = \{ v \in H^1(\Omega)|\; v_{|T} \in \mathbb{P}_l \; \forall T \in \mathbb{T}_r \} 
		\]
		the finite element space which is a subspace of $H^1(\Omega)$, based on the shape function space of piecewise polynomials of degree $l$. The shape function space is defined over the mesh $\mathbb{T}_r = \bigcup_i T_i$, where the cells $T_i$ are triangles. These spaces can be scalar-valued, vector-valued or symmetric tensor-valued, depending on the variables to be discretized. The parameter $r$ is the average size of a mesh element. In the following all the co-energy variables $\bm{e}$  are discretized by the same finite element space. The analysis were conducted using two different finite element spaces:
		\begin{enumerate}
			\item the first order Lagrange polynomials $H_r^1(\mathbb{P}_1)$;
			\item the second order Lagrange polynomials $H_r^1(\mathbb{P}_2)$.
		\end{enumerate}
		{The Lagrange multipliers $\bm{\lambda}$ are discretized by using  Lagrange polynomials defined over the boundary part where Dirichlet conditions apply. The order of the Lagrange polynomials is the same as the one chosen for the co-energy variable. The corresponding finite element space is denoted by $H_r^1(\mathbb{P}_l, \partial \Omega)$.}
		
		\subsection{Eigenvalues Computation}
		The test case for this analysis is a simple square plate of side $L$, a benchmark problem which has been studied in \cite{dawe1980rayleigh, huang1984nine, duran1999approximation} for different boundary conditions: 
		\begin{itemize}
			\item all clamped CCCC, i.e. $w_t = 0, \, \omega_n = 0,  \, \omega_s = 0$; 
			\item simply supported hard SSSS, i.e. $w_t = 0, \, M_{nn} = 0,  \, \omega_s = 0$;
			\item half-clamped half-simply supported SCSC;
			\item all clamped but one side free CCCF (for the free condition it holds $q_n = 0, \, M_{nn} = 0,  \, M_{ns} = 0$). 
		\end{itemize}
		Two different thickness to length ratios $h/L$ are considered: $h/L=0.1$ representative of a thick plate and $h/L=0.01$ representative of a thin plate. In order to compare our results with \cite{dawe1980rayleigh}, the frequencies $\omega_{mn}^h$ are computed in the following non-dimensional form:
		\begin{equation*}
		\widehat{\omega}_{mn}^h = \omega_{mn}^h L \left(\frac{2 (1 + \nu) \rho}{E}\right)^{1/2},
		\end{equation*}
		$m$ and $n$ being the numbers of half-waves occurring in the mode shapes in the $x$ and $y$ directions, respectively. The only parameters which influence the results are Poisson's ratio $\nu=0.3$, the correction factor $k$, whose value is taken equal to $k = 0.8601$ for CCCC and CCCF, $k = 0.8333$ for SSSS, $k = 0.822$ for SCSC (for comparison purposes) and the thickness-to-span ratio $h/L$. The reported non-dimensional frequencies are independent of the remaining geometrical and physical parameters.  The error is computed as:
		\begin{equation*}
		\varepsilon = \frac{\text{abs}(\widehat{\omega}_{mn}^h - \omega_{mn}^{DR})}{\omega_{mn}^{DR}},
		\end{equation*}
		where $\omega_{mn}^{DR}$ are the eigenvalues calculated in \cite{dawe1980rayleigh} using an analytical procedure. The results obtained using $H_r^1(\mathbb{P}_1), H_r^1(\mathbb{P}_2)$ for the thick and thin case are reported in Tables \ref{tab:thick}, \ref{tab:thin} respectively, together with the results from \cite{dawe1980rayleigh} (D-R) and \cite{huang1984nine} (H-H). For linear polynomials the mesh consists of a regular grid with 10 and 20 elements by side. For quadratic polynomials the mesh consists of a regular grid with 5 and 10 elements by side. Ten elements in $H_r^1(\mathbb{P}_1)$ and 5 in $H_r^1(\mathbb{P}_2)$ have the same number of degree of freedom $n=968$. Analogously, 20 elements in $H_r^1(\mathbb{P}_1)$ and 10 in $H_r^1(\mathbb{P}_2)$ have the same degree of freedom $n=3528$.  In the thick plate case,  the error is limited to $2\%$ in each for any choice of the finite element space.  On the contrary the thin case exhibits worse results, particularly for the polynomials of order 1. This must be linked to the shear locking phenomenon since $\bm{M}_{w} \propto h$, $\bm{M}_{\theta} \propto h^3$, $\bm{M}_{\kappa} \propto h^{-3}$ and $\bm{M}_{\gamma} \propto h^{-1}$. Higher order elements (like the second order Lagrange polynomials $\mathbb{P}_2$) allow alleviating the problem, but this issue is still present with very thin plate as the results worsen when the thickness to length ratio decreases to $h/L = 0.01$. Anyway, the computed eigenvalues are consistent with those obtained with other methods. The first four eigenvectors for the different cases are reported as well (see Figs. \ref{fig:CCCC} to \ref{fig:CCCF}).
	}
	
	\begin{table}[p]
		\centering	
		\begin{tabular}{|c||c||g|g||g|g||c||b|}
			\hline 
			BCs					  & Mode & \cellcolor{white}$N=10 \, (\mathbb{P}_1)$&\cellcolor{white} $N=20 \, (\mathbb{P}_1)$& \cellcolor{white}$N=5 \, (\mathbb{P}_2)$&\cellcolor{white} $N=10 \, (\mathbb{P}_2)$& H-H &\cellcolor{white} D-R \\ 
			\hline 
			\multirow{4}{*}{CCCC} &$\widehat{\omega}_{11}$& 1.5999& 1.5917 & 1.5976& 1.5914& 1.591& 1.594\\
			&$\widehat{\omega}_{21}$& 3.0615& 3.0410 & 3.0584&	3.0405& 3.039& 3.046\\
			&$\widehat{\omega}_{12}$& 3.0615& 3.0410 & 3.0677&	3.0405& 3.039& 3.046\\
			&$\widehat{\omega}_{22}$& 4.3161& 4.2682 & 4.3109&	4.2662& 4.263& 4.285\\
			\hline	
			\multirow{4}{*}{SSSS} &$\widehat{\omega}_{11}$& 0.9324& 0.9324 & 0.9304&0.9302& 0.930& 0.930\\
			&$\widehat{\omega}_{21}$& 2.2227& 2.2226& 2.2223& 2.2194& 2.219& 2.219\\
			&$\widehat{\omega}_{12}$& 2.2227& 2.2226& 2.2224& 2.2194& 2.219& 2.219\\
			&$\widehat{\omega}_{22}$& 3.4142& 3.3608& 3.4128& 3.4061& 3.405& 3.406\\
			\hline		
			\multirow{4}{*}{SCSC} &$\widehat{\omega}_{11}$& 1.3111& 1.3013& 1.3053&	1.3004 & 1.300& 1.302\\
			&$\widehat{\omega}_{21}$& 2.4155& 2.3966& 2.4040& 2.3946& 2.394& 2.398\\
			&$\widehat{\omega}_{12}$& 2.9082& 2.8871& 2.9060& 2.8858& 2.885& 2.888\\
			&$\widehat{\omega}_{22}$& 3.8906& 3.8458& 3.8721& 3.8415& 3.839& 3.852\\
			\hline
			\multirow{5}{*}{CCCF} &$\widehat{\omega}_{\frac{1}{2}1}$& 1.0855& 1.0982& 1.0845& 1.0797& 1.081&	1.089\\
			&$\widehat{\omega}_{\frac{3}{2}1}$& 1.7636& 1.7461& 1.7559& 1.7425&	1.744&	1.758\\
			&$\widehat{\omega}_{\frac{1}{2}2}$& 2.6696& 2.6575& 2.6762& 2.6547&	2.657&	2.673\\
			&$\widehat{\omega}_{\frac{5}{2}1}$& 3.2248& 3.1997& 3.2186& 3.1954&	3.197&	3.216\\
			\hline 
		\end{tabular}
		\caption{Eigenvalues for $h/L = 0.1$ using $\mathbb{P}_1$ and $\mathbb{P}_2$:	\\
			\sqbox{Blue} reference, \, \sqbox{Green} $\varepsilon<2\%$.
		}
		\label{tab:thick}
	\end{table}

	\begin{table}[p]
		\centering	
		\begin{tabular}{|c||c||c|c||m|g||c||b|}
			\hline 
			BCs					  & Mode & \cellcolor{white}$N=10 \, (\mathbb{P}_1)$&\cellcolor{white} $N=20 \, (\mathbb{P}_1)$& \cellcolor{white}$N=5 \, (\mathbb{P}_2)$&\cellcolor{white} $N=10 \, (\mathbb{P}_2)$& H-H &\cellcolor{white} D-R \\ 
			\hline 
			\multirow{4}{*}{CCCC} &$\widehat{\omega}_{11}$& \cellcolor{Yellow}0.1967&\cellcolor{Green}0.1765& 0.1872&	0.1762 &	0.1754&	0.1754\\
			&$\widehat{\omega}_{21}$& \cellcolor{Yellow}0.4030&	\cellcolor{Green}0.3604& \cellcolor{GreenYellow}0.3725&	0.3598&	0.3574&	0.3576\\
			&$\widehat{\omega}_{12}$& \cellcolor{Yellow}0.4030&	\cellcolor{Green}0.3604& 0.4055&0.3598&	0.3574&	0.3576\\
			&$\widehat{\omega}_{22}$& \cellcolor{Orange}0.6431&	\cellcolor{Green}0.5358& 0.6043&0.5335&	0.5264&	0.5274\\
			\hline 	
			\multirow{4}{*}{SSSS} &$\widehat{\omega}_{11}$&	\cellcolor{Red}0.1706&\cellcolor{Orange}0.1128& \cellcolor{Green}0.0963&	0.0963 &	0.0963&	0.0963\\
			&$\widehat{\omega}_{21}$&	\cellcolor{Purple}0.3576&\cellcolor{Yellow}0.2660& \cellcolor{Green}0.2422&	0.2406&	0.2406&	0.2406\\
			&$\widehat{\omega}_{12}$&	\cellcolor{Purple}0.3576&\cellcolor{Yellow}0.2660& \cellcolor{Green}0.2430&	0.2406&	0.2406&	0.2406\\
			&$\widehat{\omega}_{22}$&	\cellcolor{Purple}0.5803&\cellcolor{Orange}0.4442& \cellcolor{Green}0.3874&	0.3848&	0.3847&	0.3848\\
			\hline 			
			\multirow{4}{*}{SCSC} &$\widehat{\omega}_{11}$&	\cellcolor{Purple}0.1864&\cellcolor{Yellow}0.1487& 0.1492&	0.1418&	0.1411&	0.1411\\
			&$\widehat{\omega}_{21}$&	\cellcolor{Purple}0.3649&\cellcolor{Yellow}0.2829& 0.2827&	0.2683&	0.2668&	0.2668\\
			&$\widehat{\omega}_{12}$&	\cellcolor{Orange}0.3987&\cellcolor{GreenYellow}0.3485& 0.3608&	0.3394&	0.3377&	0.3377\\
			&$\widehat{\omega}_{22}$&	\cellcolor{Purple}0.6075&\cellcolor{Yellow}0.4933& 0.4940&	0.4654&	0.4604&	0.4608\\
			\hline 		
			\multirow{5}{*}{CCCF} &$\widehat{\omega}_{\frac{1}{2}1}$& \cellcolor{Yellow}0.1238& \cellcolor{Green}0.1166& \cellcolor{GreenYellow}0.1197& 0.1169&0.1166&0.1171\\
			&$\widehat{\omega}_{\frac{3}{2}1}$& \cellcolor{Yellow}0.2207& \cellcolor{Green}0.1954& 0.2092& 0.1960&	0.1949&	0.1951\\
			&$\widehat{\omega}_{\frac{1}{2}2}$& \cellcolor{GreenYellow}0.3204& \cellcolor{Green}0.3078& \cellcolor{GreenYellow}0.3188& 0.3089&	0.3080&	0.3093\\
			&$\widehat{\omega}_{\frac{5}{2}1}$& \cellcolor{Yellow}0.4144& \cellcolor{Green}0.3751& 0.3938& 0.3757&	0.3736&	0.3740\\			
			\hline 
		\end{tabular} 	
		\caption{Eigenvalues for $h/L = 0.01$ using $\mathbb{P}_1$ and $\mathbb{P}_2$: \\
			\sqbox{Blue} reference, \, \sqbox{Green} $\varepsilon<2\%$,\, \sqbox{GreenYellow} $\varepsilon<5\%$, \, 	
			\sqbox{Yellow} $\varepsilon<15\%$, \, \sqbox{Orange} $\varepsilon<30\%$, \, \sqbox{Purple} $\varepsilon<50\%$,	\, \sqbox{Red} $\varepsilon<80\%$. }
		\label{tab:thin}
	\end{table}

\subsection{Time-domain Simulations}
\begin{figure}[t]%
	\centering
	\subfloat[][$w(t = 0.25 \, t_{\text{fin}})$]{%
		\label{fig:sim1-a}%
		\includegraphics[width=0.45\textwidth]{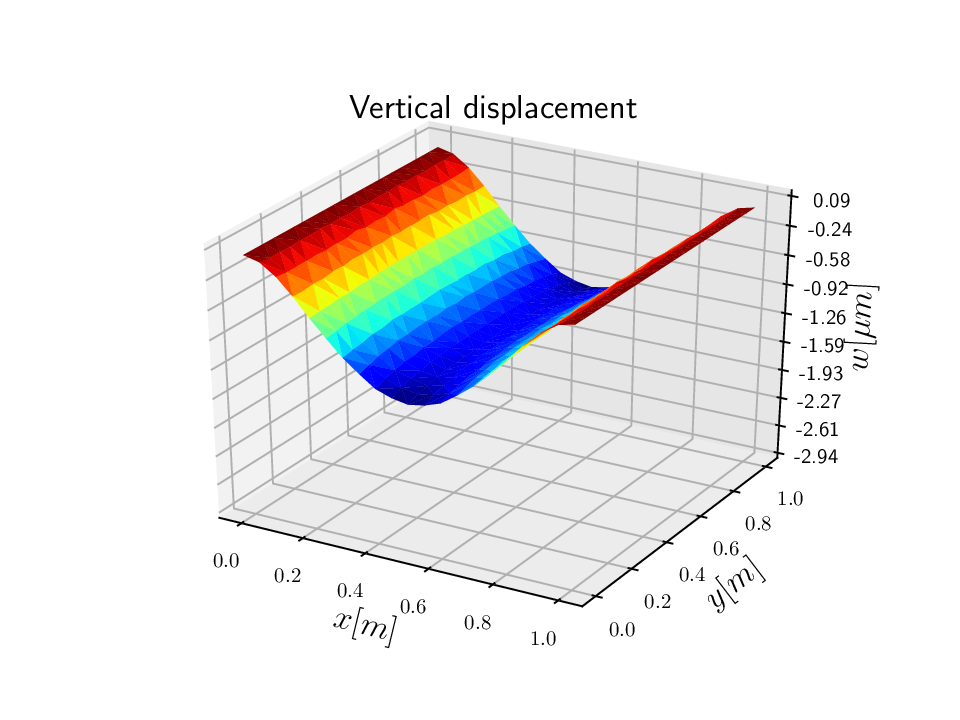}}%
	\hspace{8pt}%
	\subfloat[][$w(t = 0.50 \, t_{\text{fin}})$]{%
		\label{fig:sim1-b}%
		\includegraphics[width=0.45\textwidth]{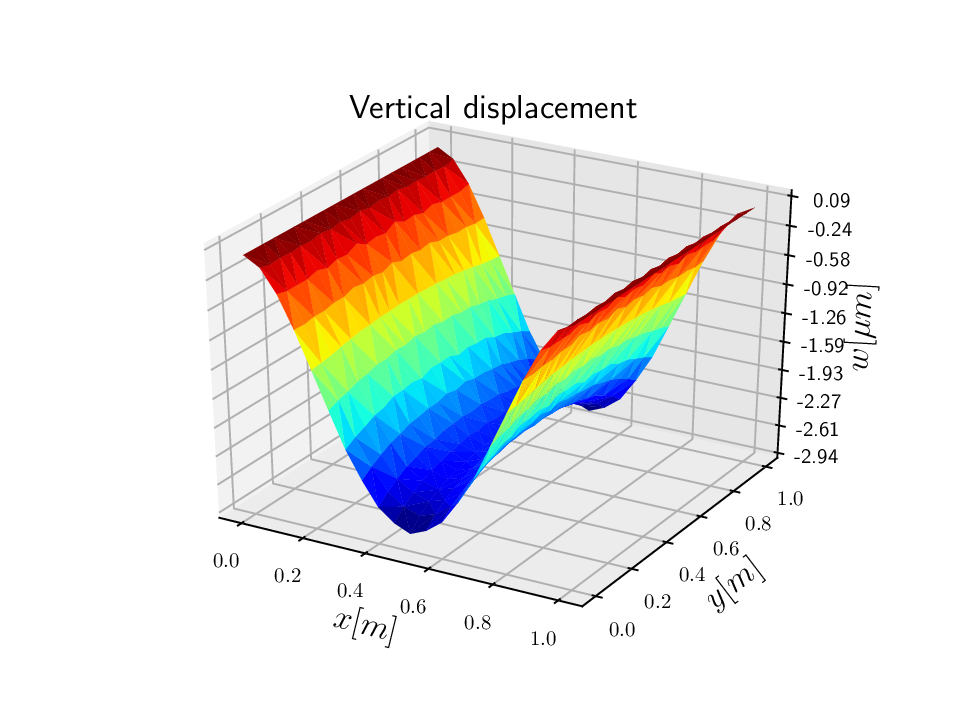}} \\
	\subfloat[][$w(t = 0.75 \, t_{\text{fin}})$]{%
		\label{fig:sim1-c}%
		\includegraphics[width=0.45\textwidth]{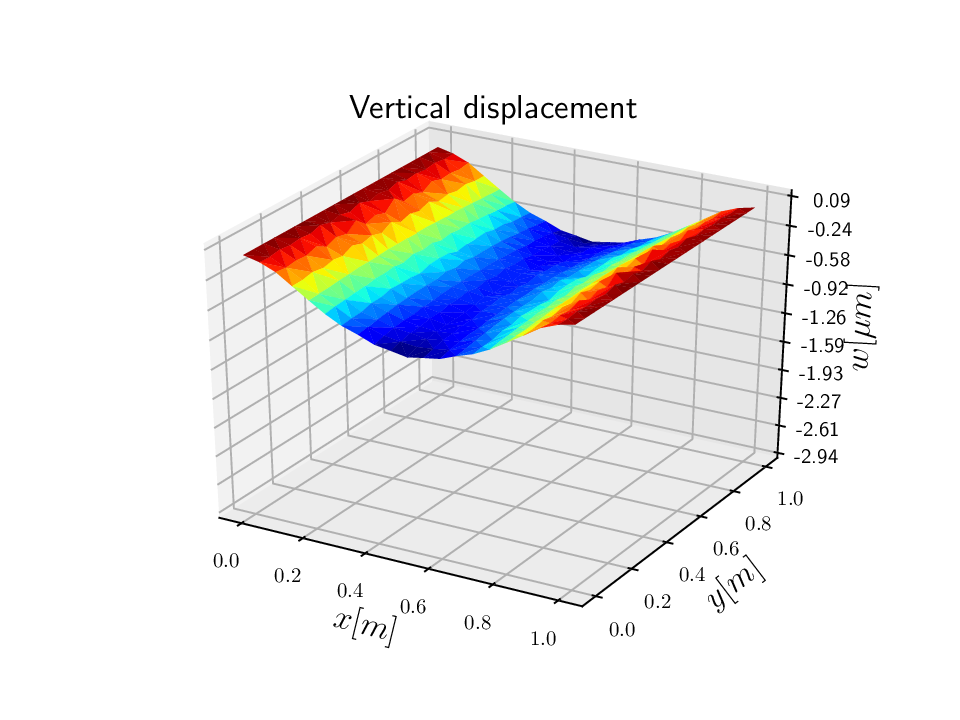}}%
	\hspace{8pt}%
	\subfloat[][$w(t = t_{\text{fin}})$]{%
		\label{fig:sim1-d}%
		\includegraphics[width=0.45\textwidth]{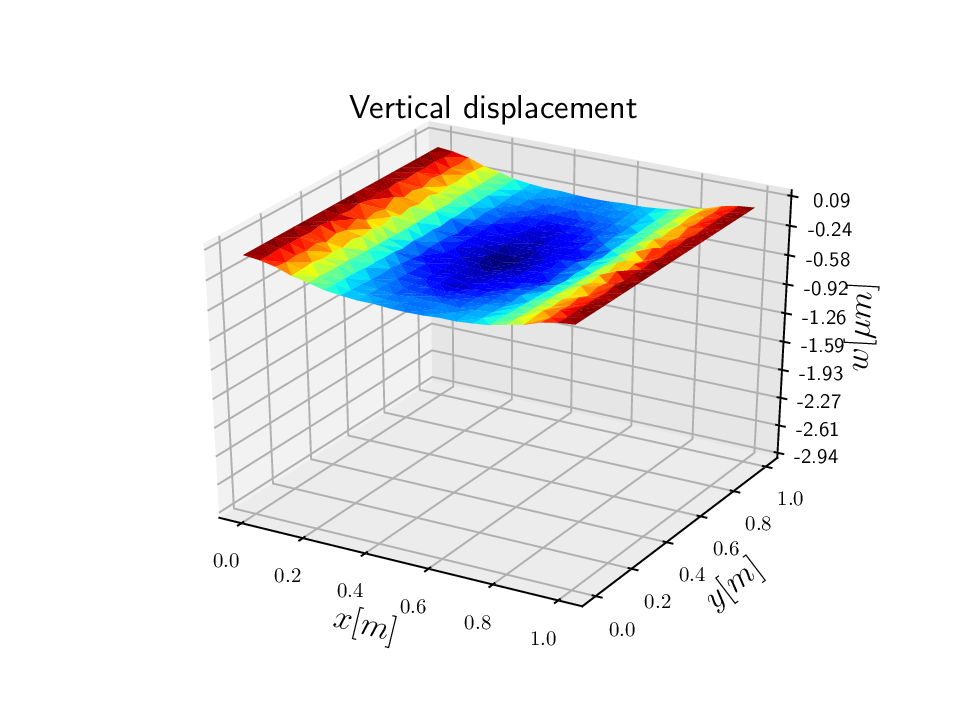}}%
	\caption[Snapshots of the displacement field]{Snapshots for Simulation $n^\circ 1$.}%
	\label{fig:sim1}%
\end{figure}

\begin{figure}[h]%
	\centering
	\subfloat[][$w(t = 0.25 \, t_{\text{fin}})$]{%
		\label{fig:sim2-a}%
		\includegraphics[width=0.45\textwidth]{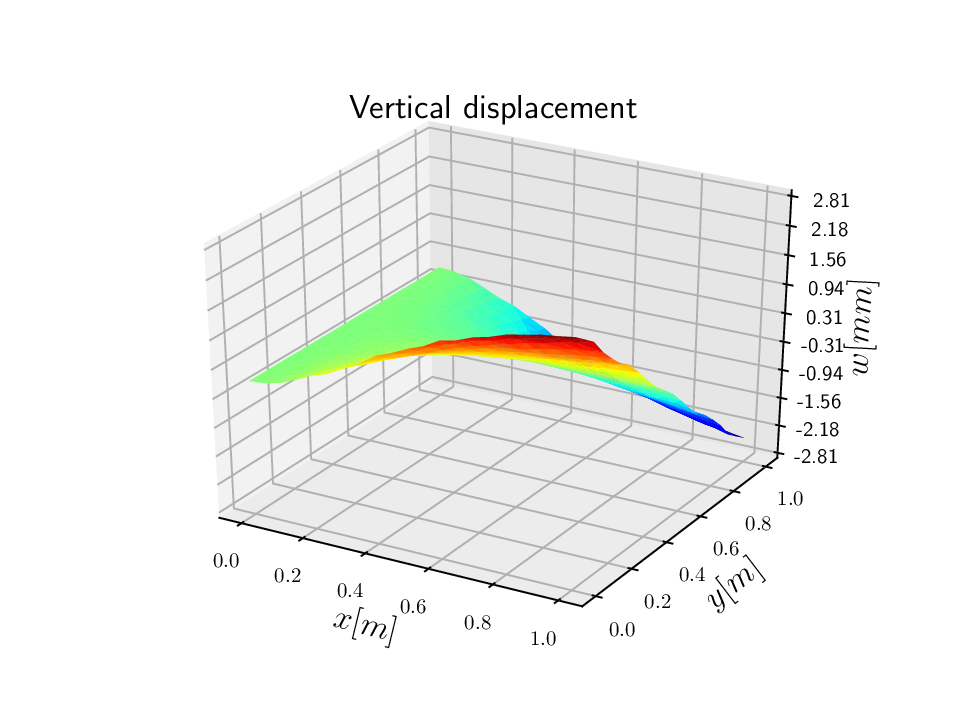}}%
	\hspace{8pt}%
	\subfloat[][$w(t = 0.50 \, t_{\text{fin}})$]{%
		\label{fig:sim2-b}%
		\includegraphics[width=0.45\textwidth]{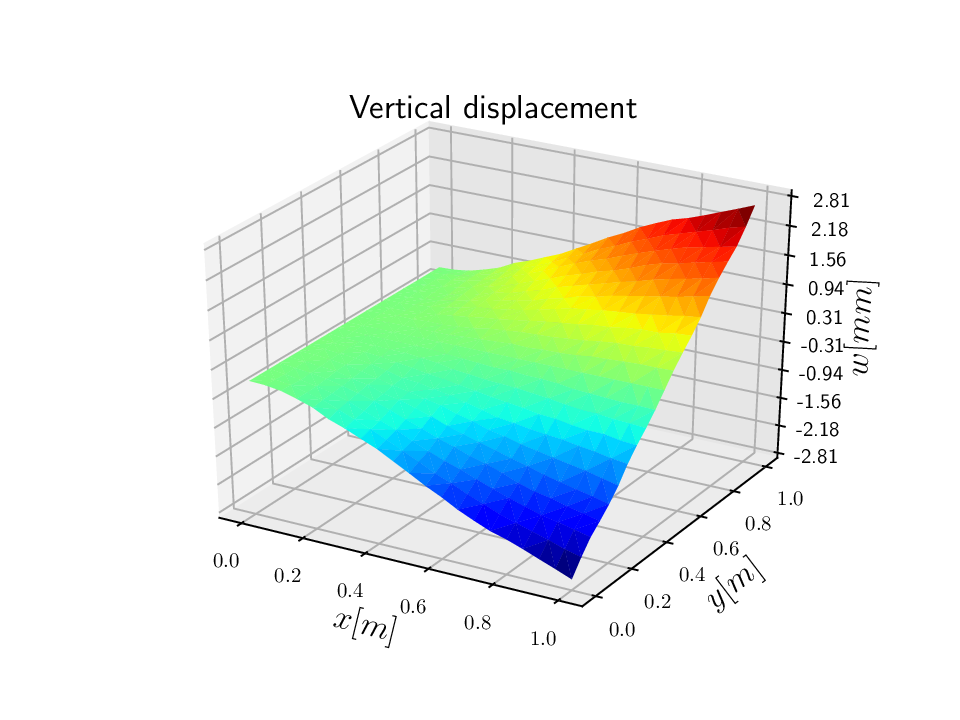}} \\
	\subfloat[][$w(t = 0.75 \, t_{\text{fin}})$]{%
		\label{fig:sim2-c}%
		\includegraphics[width=0.45\textwidth]{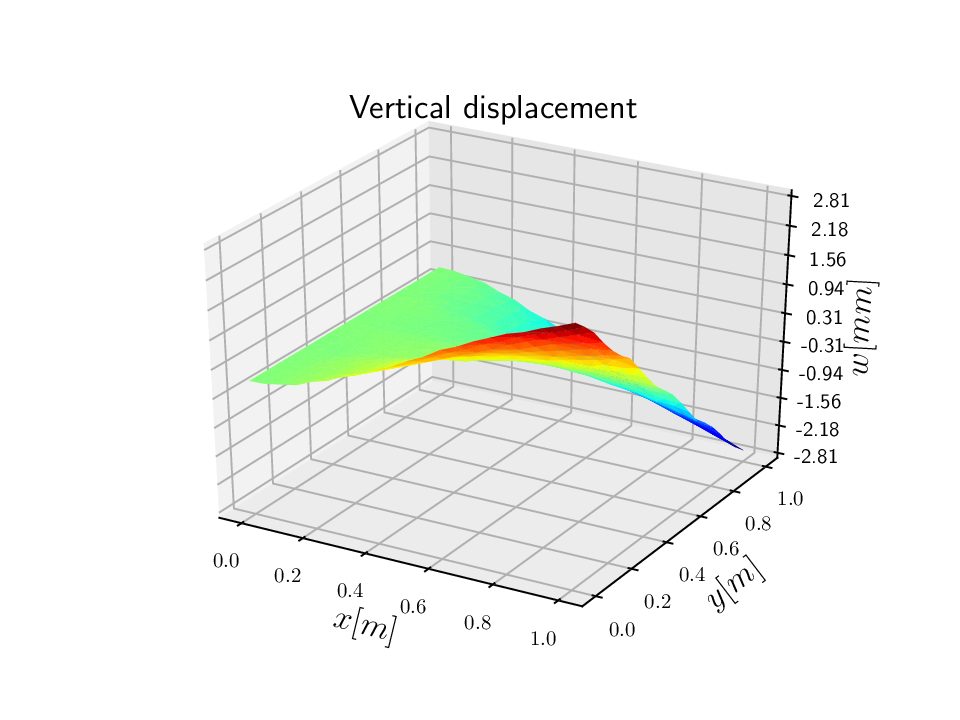}}%
	\hspace{8pt}%
	\subfloat[][$w(t = t_{\text{fin}})$]{%
		\label{fig:sim2-d}%
		\includegraphics[width=0.45\textwidth]{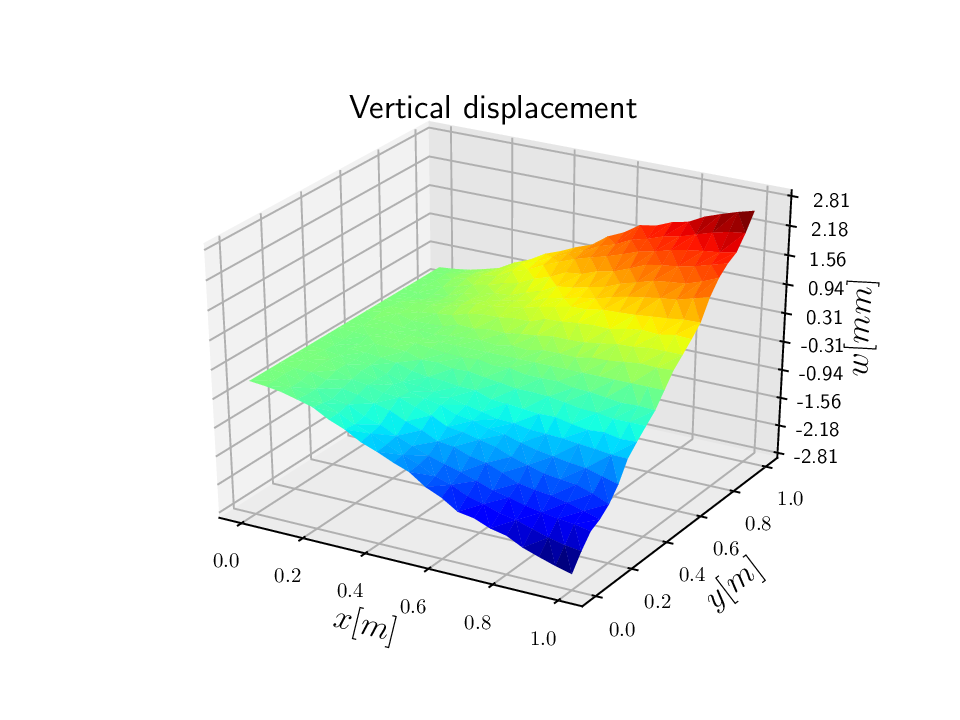}}%
	\caption[Snapshots of the displacement field]{Snapshots for Simulation $n^\circ 2$.}%
	\label{fig:sim2}%
\end{figure}

\begin{table}[h!]
	\centering
	\begin{tabular}{|c|c|}
		\hline 
		\multicolumn{2}{|c|}{Plate Parameters} \\ 
		\hline 
		$E$ & 70 $[GPa]$ \\ 
		$\rho$ & $2700\; [kg/m^3]$ \\ 
		$\nu$& 0.35 \\ 
		$k$& 5/6 \\ 
		$h/L$& 0.1 \\ 
		$L$& 1 $[m]$ \\
		\hline 
	\end{tabular} 
	\begin{tabular}{|c|c|}
		\hline 
		\multicolumn{2}{|c|}{Simulation} \\  
		\multicolumn{2}{|c|}{Settings} \\  
		\hline 
		Integrator & St\"ormer-Verlet \\ 
		$\Delta t$ & $0.001\; [ms]$ \\ 
		$t_{\text{end}}$ & $10\; [ms]$ \\ 
		N$^\circ$ Elements & 10 \\
		FE space & $H_{r = L/10}^1(\mathbb{P}_2)\text{ for }\bm{e} \times H_{r = L/10}^1(\mathbb{P}_2, \partial\Omega)\text{ for }\bm{\lambda}$ \\
		\hline 
	\end{tabular} 
	\captionsetup{width=0.95\linewidth}
	\vspace{1mm}
	\captionof{table}{Physical parameters and simulations settings.}
	\label{tab:par}
\end{table}
	\begin{figure}[t]%
		\centering
		\subfloat[][Simulation $n^\circ 1$]{%
			\label{fig:sim1-H}%
			\includegraphics[width=0.47\columnwidth]{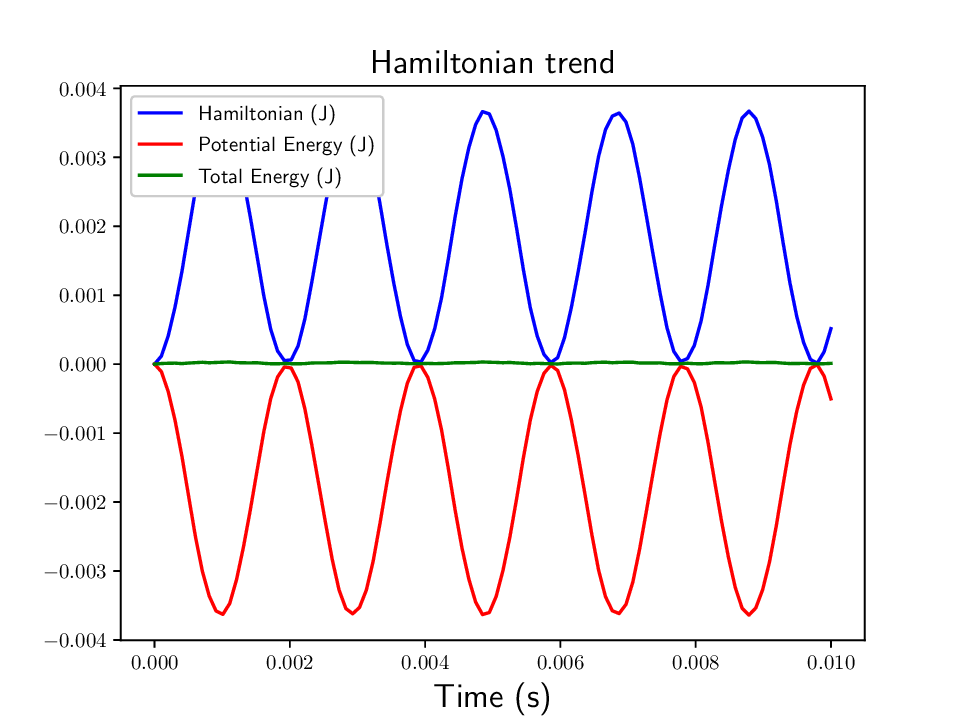}}%
		\hspace{8pt}%
		\subfloat[][Simulation $n^\circ 2$]{%
			\label{fig:sim2-H}%
			\includegraphics[width=0.45\columnwidth]{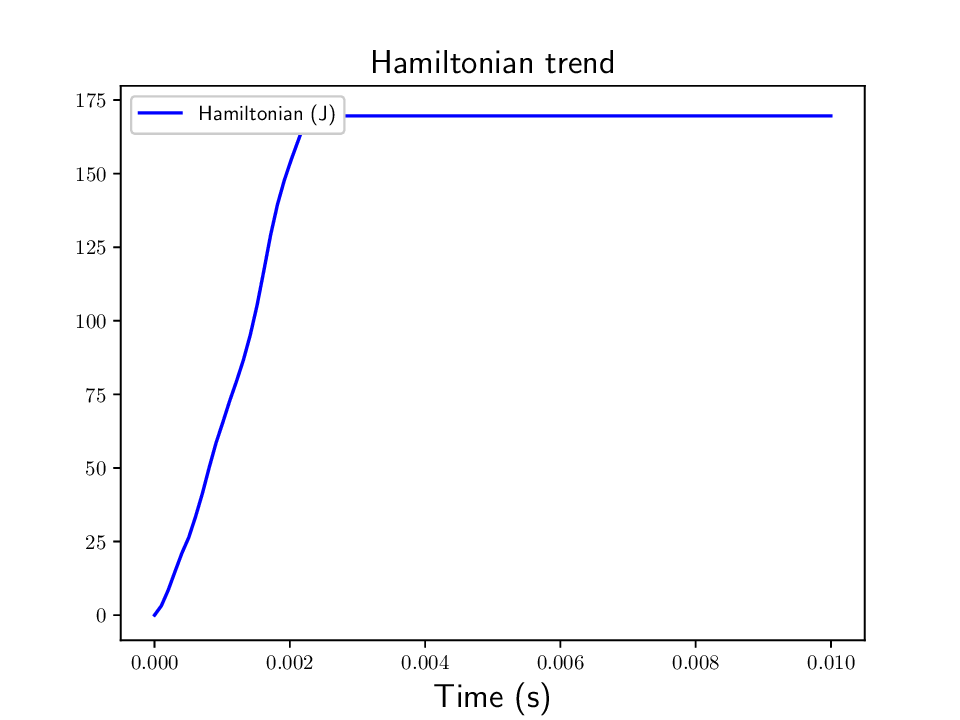}} \\
		\caption[Hamiltonian]{Hamiltonian trend for the two simulations.}
		\label{fig:Hamiltonian}
	\end{figure}
In this analysis, a square plate, subject either to a non null shear force on the boundaries, either to a distributed force over the domain is considered. The physical parameters and simulation settings are reported in Table \ref{tab:par}. Ten elements for each side and second order Lagrange polynomials are used. The St\"ormer-Verlet time integrator is employed, {so that the symplectic structure is preserved.} Two different simulations with different boundary conditions are considered. The initial conditions are set to zero for each variable. For the first simulation, a cantilever plate subject to gravity is considered. In order to add gravity the corresponding functional has to be discretized:
\[\bm{f}_{\text{gravity}} = - \int_{\Omega} v_w \rho h g \d{\Omega},
\]
where $g = 10 [m/s^2]$ is the gravity acceleration. This excitation, that applies from $t=0$, is then simply added to the first line of system \eqref{eq:PHdiscr_mixed}.
For the first simulation the following boundary conditions, corresponding to the CFCF cases, are considered:
\begin{equation*} \text{Simulation $n^\circ \, 1$} \;
\begin{cases}
w_t = 0, \, \omega_n = \omega_s= 0, \quad &\text{for } x=0 \text{ and } x=1,\\
q_n = 0, \, M_{nn} = M_{ns}= 0, \quad &\text{for } y=0,\\
q_n = 0, \, M_{nn} = M_{ns}= 0, \quad &\text{for } y=1.\\
\end{cases}
\end{equation*}
Since this force admits a potential, the Hamiltonian does not correspond to the the total energy, that now includes the potential energy:
\begin{equation*}
E_p = \int_{\Omega} \rho h g w \d{\Omega}, 
\end{equation*}
where $w$ is the vertical displacement field. \\
For the second simulation, the following  boundary conditions are considered:
\begin{equation*} \text{Simulation $n^\circ \, 2$} \;
\begin{cases}
w_t = 0, \, \omega_n = \omega_s= 0, \quad &\text{for } x=0,\\
q_n = 0, \, M_{nn} = M_{ns}= 0, \quad &\text{for } x=1,\\
q_n = +f(x, t), \, M_{nn} = M_{ns}= 0, \quad &\text{for } y=0,\\
q_n = -f(x, t), \, M_{nn} = M_{ns}= 0, \quad &\text{for } y=1,\\
\end{cases}
\end{equation*}
where $f(x, t)$ is the excitation defined by:
\begin{equation*}
f(x, t) = \begin{cases}
10^5 \,\sin \left( \frac{\pi}{L} x \right) \; [Pa \cdot m], \; &\forall t< 0.25 \, t_{\text{end}}, \\
0, \; &\forall t\geq 0.25 \, t_{\text{end}}, \\
\end{cases}
\end{equation*} 
$t_{\text{end}}$ being defined in Table \ref{tab:par}. In this case inhomogeneous boundary conditions are considered. \\
Snapshots of the vertical displacement field are reported in Figures \ref{fig:sim1}, \ref{fig:sim2}. This field is obtained from the velocity field $e_w = \diffp{w}{t}$ by applying the trapezoidal rule integration. For both simulations the output is consistent with the imposed BC and with the physical intuition of the observed phenomenon. The symplectic integration has been used to demonstrate numerically the conservation of total energy. The Hamiltonian for both simulations is reported in Fig. \ref{fig:Hamiltonian}. For the first one the total energy $E_t = H + E_p$ is the conserved quantity and it remains constant equal to zero as expected. For the second simulation the Hamiltonian is the conserved quantity, once the excitation is set back to zero. The conservation of energy is also obvious from Fig. \ref{fig:sim2}. The boundary conditions are such that the energy increases monotonically as long as the excitation is present.  The deformation attained at the final time of the excitation  repeats itself symmetrically with respect to axis $y = L/2$ with period given by $1/2 \; t_{\text{end}}$.

\section*{Conclusions and Perspectives}

In this paper the port-Hamiltonian formulation of the Mindlin plate was enriched by the equivalent tensorial formulation and by a symplectic discretization method. The PFEM method, applied to this formulation, proves again its versatility, since it stays valid and applicable even in more complicated examples like the one presented in this article. Many features of this method and of the tensorial PH formulation are of interest:
\begin{itemize}
	\item its capability of preserving the port-Hamiltonian structure;
	\item the natural derivation of boundary port-variables as inputs;
	\item the possibility of dealing with mixed boundary conditions inside the framework of  port-Hamiltonian descriptor systems PHDAEs,  detailed in \cite{beattie2018linear};
	\item the easy implementability of the method using standard Finite Element libraries (Fenics \cite{LoggMardalEtAl2012} in our case);
	\item a coordinate free representation which makes possible to treat other planar geometries (circular, ellipsoidal, annular etc).
\end{itemize}  
The formulation presented in this paper could be completed with a precise analysis of the well-posedness, in the input-output sense (as in \cite{waveEqZwart} for the wave equation in $\mathbb{R}^d$). Furthermore, a proper convergence analysis is needed. {The Arnold-Winther element \cite{Arnold2002} should be investigated as they provide a conforming approximation of space $H^{\text{Div}}(\Omega, \mathbb{R}^{2 \times 2}_\text{sym})$. Unfortunately, this are not included inside FEniCS (or in any standard library)}.  \\ \\
Another important aspect is the implementation of suitable control laws. The proposed method paves the way to the use of passivity-based approaches and of energy shaping methods. These techniques have already been largely studied in the literature \cite{Ortega2002, OrtegaContrInt} for the finite-dimensional case. However, the distributed case is mainly limited to one geometrical dimension \cite{MaccContrDist}.  It would be of great interest, especially for control engineers, to conceive a controller able to respect given performance specifications. The port-Hamiltonian formalism and its powerfulness in modeling complex systems could be linked to standard control methodologies, already well known in the industrial field, like the $LQR$ or $H^\infty$ methodologies.   

\section*{Acknowledgments}
This work is  supported by the project ANR-16-CE92-0028,
entitled {\em Interconnected Infinite-Dimensional systems for Heterogeneous
	Media}, INFIDHEM, financed by the French National
Research Agency (ANR) and the Deutsche Forschungsgemeinschaft (DFG).
Further information is available at \\ {\url{https://websites.isae-supaero.fr/infidhem/the-project}}.  \\
Moreover the authors would like to thank Michel Sala\"un from ISAE-SUPAERO for the fruitful and insightful discussions.

	\bibliographystyle{unsrt}
	\bibliography{biblio_Mindlin_Revision} 
	
\end{document}